\documentclass[11pt]{article}
\usepackage{amsfonts,mathrsfs}
\usepackage{amssymb}
\usepackage{amsmath}
\usepackage{hyperref}
\usepackage{color}
=0pt
\def\sqr#1#2{{\vcenter{\vbox{\hrule height.#2pt
              \hbox{\vrule width.#2pt height#1pt \kern#1pt \vrule
width.#2pt}
              \hrule height.#2pt}}}}
\def\signed #1{{\unskip\nobreak\hfil\penalty50
              \hskip2em\hbox{}\nobreak\hfil#1
              \parfillskip=0pt \finalhyphendemerits=0 \par}}
\def\endpf{\signed {$\sqr69$}}
\def\dbE{{\mathbb{E}}}
\def\dbF{{\mathbb{F}}}

\def\dbP{{\mathbb{P}}}

\def\dbR{{\mathbb{R}}}
\def\dbS{{\mathbb{S}}}

\def\dbW{{\mathbb{W}}}

\def\d{\delta}
\def\e{\varepsilon}

\def\k{\kappa}

\def\si{\sigma}

\def\f{\varphi}

\def\om{\omega}

\def\3n{\negthinspace \negthinspace \negthinspace }
\def\2n{\negthinspace \negthinspace }
\def\1n{\negthinspace }
\def\ns{\noalign{\smallskip} }

\def\ds{\displaystyle}
\def\D{\Delta}

\def\Om{\Omega}
\def\om{\omega}
\def\cC{{\cal C}}

\def\cF{{\cal F}}

\def\cL{{\cal L}}

\def\cO{{\cal O}}

\def\cU{{\cal U}}

\def\cX{{\cal X}}
\def\cY{{\cal Y}}

\def\mE{{\mathbb{E}}}

\def\no{\noindent}

\def\ss{\smallskip}
\def\ms{\medskip}
\def\bs{\bigskip}
\def\q{\quad}
\def\qq{\qquad}
\def\hb{\hbox}

\def\lan{\big\langle}
\def\ran{\big\rangle}

\def\pa{\partial}

\def\wt{\widetilde}
\def\cd{\cdot}

\def\ae{\hbox{\rm a.e.{ }}}

\def\deq{\mathop{\buildrel\D\over=}}

\def\({\Big (}
\def\){\Big )}
\def\[{\Big[}
\def\]{\Big]}

\def\={\buildrel \triangle \over =}

\def\be{\begin{equation}}
\def\bel{\begin{equation}\label}
\def\ee{\end{equation}}
\def\bea{\begin{eqnarray}}
\def\eea{\end{eqnarray}}
\def\bt{\begin{theorem}}
\def\et{\end{theorem}}
\def\bc{\begin{corollary}}
\def\ec{\end{corollary}}
\def\bl{\begin{lemma}}
\def\el{\end{lemma}}
\def\bp{\begin{proposition}}
\def\ep{\end{proposition}}
\def\br{\begin{remark}}
\def\er{\end{remark}}
\def\ba{\begin{array}}
\def\ea{\end{array}}
\def\bde{\begin{definition}}
\def\ede{\end{definition}}
\def\bex{\begin{example}}
\def\eex{\end{example}}

\newtheorem{lemma}{Lemma}[section]
\newtheorem{remark}{Remark}[section]
\newtheorem{example}{Example}[section]
\newtheorem{theorem}{Theorem}[section]
\newtheorem{corollary}{Corollary}[section]

\newtheorem{definition}{Definition}[section]
\newtheorem{proposition}{Proposition}[section]

\DeclareMathOperator*\uplim{\overline{lim}}

\allowdisplaybreaks
 \makeatletter
   
   \@addtoreset{equation}{section}
\makeatother
\begin{document}

\title{\bf Stochastic Verification Theorem for
Infinite Dimensional Stochastic Control Systems}

\author{Liangying Chen\footnote{School
of Mathematics, Sichuan University, Chengdu, P.
R. China, and Sorbonne Universit\'es,  UPMC Univ
Paris 06, Paris, France. Email:
chenli@ljll.math.upmc.fr. Liangying Chen is
supported by the European Union's Horizon 2020
research and innovation programme under the
Marie Sklodowska-Curie grant agreement No
945322.} ~~~ and ~~~ Qi L\"{u}\footnote{School
of Mathematics, Sichuan University, Chengdu, P.
R. China. Email: lu@scu.edu.cn. Qi L\"u is
supported by the NSF of China under grants
12025105 and 11971334.}}
\date{}

\maketitle

\begin{abstract}
The verification theorem serving as an optimality condition for the
optimal control problem, has been expected and studied for a long
time. The purpose of this paper is to establish this theorem for
control systems governed by stochastic evolution equations in
infinite dimensions, in which both the drift and the diffusion terms
depend on the controls.
\end{abstract}

\bs

\no{\bf 2010 Mathematics Subject
Classification}. 93E20.

\bs

\no{\bf Key Words:} Optimal control, value function, stochastic distributed
parameter systems, stochastic verification theorem

\bs

\section{Introduction}

We begin with some notations. Let $T>0$, and let $(\Omega,
\mathcal{F}, \mathbb{P})$ be a complete probability space, on which
a separable Hilbert space $\wt H$-valued cylindrical Brownian motion
$W(\cdot)$ is defined. Denote by $\textbf{F}\deq \{\cF_t\}_{t\geq
0}$ the natural filtration generated by $W(\cdot)$ and by
$\mathbb{F}$ the progressive $\sigma$-algebra with respect to
$\textbf{F}$.

Let $\cX$ be a Banach space. For any $t\in[0,T]$ and $p\in
[1,\infty)$, denote by $L_{\cF_t}^p(\Om;\cX)$ the Banach space of
all $\cF_t$-measurable random variables $\xi:\Om\to \cX$ such that
$\mathbb{E}|\xi|_\cX^p < \infty$, with the canonical  norm. Denote
by $L^{p}_{\dbF}(\Om;C([t,T];\cX))$ the Banach space of all
$\cX$-valued $\mathbf{F}$-adapted continuous processes
$\phi(\cdot)$, with the norm\vspace{-3mm}
$$
|\phi(\cd)|_{L^{p}_{\dbF}(\Om;C([t,T];\cX))} \=
\Big[\mE\sup_{\tau\in
    [t,T]}|\phi(\tau)|_\cX^p\Big]^{1/p}.  $$
Also, denote by $C_{\dbF}([t,T];L^{p}(\Om;\cX))$ the Banach space of
all $\cX$-valued $\mathbf{F}$-adapted processes $\phi(\cdot)$ such
that $\phi(\cdot):[t,T] \to L^{p}_{\cF_T}(\Om;\cX)$ is continuous,
with the norm
$$
|\phi(\cd)|_{C_{\dbF}([t,T];L^{p}(\Om;\cX))} \= \sup_{\tau\in
    [t,T]}\left[\mE|\phi(\tau)|_\cX^p\right]^{1/p}.
$$

Fix any $p_1,p_2 \in[1,\infty]$. Put
$$
\begin{array}{ll}
     \ds
    L^{p_2}_\dbF(t,T;L^{p_1}(\Om;\cX)) =\Big\{\f:(t,T)\times\Om\to
    \cX\;\Big|\;\f(\cd)\hb{ is $\mathbf{F}$-adapted and
    }\int_t^T\Big(\dbE|\f(\tau)|_X^{p_1}\Big)^{\frac{p_2}
        {p_1}}d\tau<\infty\Big\}.
\end{array}
$$
Clearly, $L^{p_2}_\dbF(t,T;L^{p_1}(\Om;\cX))$ is a Banach space with
the canonical norm. If $p_1=p_2$, we simply write the above spaces
as $L^{p_1}_\dbF(t,T;\cX)$. Put
$$
\begin{array}{ll}
\ds L^{p_2}_{S,\dbF}(t,T;L^{p_1}(\Om;\cX))
=\Big\{\f:(t,T)\times\Om\to \cX\;\Big|\;|\f(\cd)|_{\cX}\in
L^{p_2}_\dbF(t,T;L^{p_1}(\Om;\dbR))\Big\}.
\end{array}
$$
Similarly, if $p_1=p_2$, we simply write the above spaces as
$L^{p_1}_{S,\dbF}(t,T;\cX)$.

For $r\in[0,T]$ and $f\in L^1_{\cF_T}(\Om;\cX)$, denote by
$\dbE(f|\cF_r)$ the conditional expectation of $f$ with respect to
$\cF_r$ and by $\dbE f$ the mathematical expectation of $f$.

\ss Let $\cY$ be another Banach space. Denote by $\cL(\cX; \cY)$, or
$\cL(\cX)$ if $\cY=\cX$, the Banach space of all bounded linear
operators from $\cX$ to $\cY$ with the usual operator norm. When
$\cX$ is a Hilbert space, write $\dbS(\cX)$ for the space of all
bounded linear self-adjoint operators on $\cX$.

For $v\in C([0,T]\times \cX)$ and $(t,\eta)\in [0,T)\times \cX$, the
second-order parabolic superdifferential of $v$ at $(t,\eta)$ is
defined as follows:
\begin{eqnarray*}
D_{t,x}^{1,2,+}v(t,\eta) \3n& \deq\3n & \Big\{(r,p,P)\in
\mathbb{R}\times H\times \dbS(\cX)\Big|
\displaystyle\uplim\limits_{\begin{subarray}{1}
s\downarrow t, s\in [0,T)\\
y\to \eta
\end{subarray}}\frac{1}{|s-t|+|\eta-y|_\cX^2}\\
& & \Big[v(s,y)-v(t,\eta)-r(s-t)-\langle
p,y-\eta\rangle_\cX-\frac{1}{2}\langle
P(y-\eta),y-\eta\rangle_\cX\Big]\leq 0 \Big\}.
\end{eqnarray*}
Note that the limit in $t$ is from the right.
This fits the general irreversibility of
evolution equations.

\ss

Now we can introduce the control problem. Let $H$  be a separable
Hilbert space, and $A:D(A)\subset H\to H$ be a linear operator,
which generates a $C_{0}$-semigroup $\{S(t)\}_{t\ge 0}$ on $H$.
Write $ \mathcal{L}_2^0$ for the space of all Hilbert-Schmidt
operators from $\wt H$ to $H$, which is also a separable Hilbert
space.  Let $U$ be a separable metric space with a metric ${\bf
d}(\cdot,\cdot)$. For $t\in [0,T)$, put
\begin{equation*}
    \mathcal{U}[t,T]\deq\big\{u:[t,T]\times\Omega\to U\big| u \mbox{ is
        $\textbf{F}$-adapted}\big\}.
\end{equation*}

The control system under consideration in this paper is given as follows:
\begin{equation}\label{system1}
    \begin{cases}\ds
        dX(t)=\big(AX(t)+a(t,X(t),u(t))\big)dt+b(t,X(t),u(t))dW(t), &t \in (0,T],\\
        \ns\ds X(0)=\eta \in H,
    \end{cases}
\end{equation}
and the cost functional is
\begin{equation}\label{cost1}
    \mathcal{J}(\eta;u(\cdot))=\mathbb{E}\Big(\int_0^T
    f(t,X(t),u(t))dt+h(X(T))\Big).
\end{equation}

We make the following assumptions for the control system
\eqref{system1} and the cost functional \eqref{cost1}:

\ss

{\bf (S1)} {\it  Suppose that: i) $a(\cdot, \cdot,\cdot): [0,T]\times H\times U \to H$ is
    $\mathcal{B}([0,T])\otimes \mathcal{B}(H)\otimes
    \mathcal{B}(U)/\mathcal{B}(H)$-measurable and
    $b(\cdot,\cdot,\cdot): [0,T]\times H\times U \to
    \mathcal{L}_2^0$ is $\mathcal{B}([0,T])\otimes \mathcal{B}(H)\otimes
    \mathcal{B}(U)/\mathcal{B}(\mathcal{L}_2^0)$-measurable;  ii) for any
    $ \eta\in  H$, the maps $a(\cdot,\eta,\cdot):[0,T]\times U\to H $ and
    $b(t,\eta,\cdot): U\to \mathcal{L}_2^0$ are continuous; and iii)
    for any $(t,\eta_1,\eta_2,u)\in  [0,T]\times H \times H\times U$,
    \begin{equation*}
        \begin{cases}
            |a(t,\eta_1,u)-a(t,\eta_2,u)|_H\leq C|\eta_1-\eta_2|_H,\\
            |b(t,\eta_1,u)-b (t,\eta_2,u)|_{\mathcal{L}_2^0}\leq C|\eta_1-\eta_2|_H,\\
            |a(t,0,u)|_H \leq C, \q \ \ \ \ |b(t,0,u)|_{\mathcal{L}_2^0}\leq
            C.
        \end{cases}
\end{equation*}}

\ss

{\bf (S2)} {\it  Suppose that: i) $f(\cdot,\cdot,\cdot):[0,T]\times H\times U\to \mathbb{R}$
    is $\mathcal{B}([0,T])\otimes \mathcal{B}(H)\otimes
    \mathcal{B}(U)/\mathcal{B}(\mathbb{R})$-measurable and
    $h(\cdot):H\to \mathbb{R}$ is
    $\mathcal{B}(H)/\mathcal{B}(\mathbb{R})$-measurable; ii) For any
    $\ \eta\in  H$, the functional $f(\cdot,\eta,\cdot): [0,T]\times U\to
    \mathbb{R}$ is continuous; and iii) For any $(t,\eta_1,\eta_2,u)\in
    [0,T]\times H\times H\times U$,
    \begin{equation*}
        \begin{cases}
            |f(t,\eta_1,u)-f(t,\eta_2,u)|\leq C|\eta_1-\eta_2|_H, \\
            |h(\eta_1)-h(\eta_2)|\leq C|\eta_1-\eta_2|_H\\
            |f(t,0,u)|\leq C,  \q \ \ |h(0)|\leq C.
        \end{cases}
\end{equation*}}

Hereafter, we use $C$ to denote the generic
constant, which may change from line to line.

\begin{remark}\label{rem-3.1}
The boundedness condition on $f$ and $h$ in (S2)
is just for the convenience of computation and
to emphasize the main arguments, which can be
relaxed in the following manner:
$$
\begin{array}{ll}\ds
\big|f(t,\eta_1,u)-f(t,\eta_2,u)\big| +
\big|h(\eta_1)-h(\eta_2)\big|\leq
C\big(1+|\eta_1|_H +|\eta_2|_H\big)|\eta_1-\eta_2|_H,\\
\ns\ds\hspace{6.5cm} \forall
(t,\eta_1,\eta_2,u)\in [0,T]\times H\times
H\times U.
\end{array}
$$

\end{remark}

Under {\bf (S1)}, for any $u(\cdot)\in \mathcal{U}[0,T]$,   the control system
\eqref{system1} has a unique mild solution $X(\cdot) \in C_\dbF([0,T];$ $L^2(\Omega;H))$ (see \cite[Theorem 3.14]{Lu2021} for example).

\ss

Consider the following optimal control problem:

\ss

\textbf{Problem} $\boldsymbol{(S_{\eta})}$.  For any given $\eta\in
H$, find a
$\bar{u}(\cdot)\in \mathcal{U}[0,T]$ such that
\begin{equation}\label{OP1}
    \mathcal{J}(\bar{u}(\cdot))=\inf\limits_{u(\cdot)\in \mathcal{U}[0,T]}\mathcal{J}(u(\cdot)).
\end{equation}

Any $\bar{u}(\cdot)\in \mathcal{U}[0,T]$ satisfying \eqref{OP1} is called an {\it optimal control} (of
\textbf{Problem} $\boldsymbol{(S_{\eta})}$). The corresponding state
$\overline{X}(\cdot)$ is called an {\it optimal state}, and
$(\overline{X}(\cdot),\bar{u}(\cdot))$ is called an {\it optimal pair}.

Let us recall the stochastic dynamic programming principle for solving \textbf{Problem} $\boldsymbol{(S_{\eta})}$. In the literature, the stochastic dynamic programming principle for
\textbf{Problem} $\boldsymbol{(S_{\eta})}$ in weak formulation is
already established. A nice treatise for that is \cite{Fabbri2017}.
In that formulation, probability spaces and Brownian motions vary
with the controls. In other words, the probability space and
Brownian motion are  part of the control. Usually, optimal control
problems for SEEs are formulated in strong formulation, i.e., the
probability space and the Brownian motion are fixed. Hence, it is
natural to ask whether the stochastic dynamic programming principle
holds in strong formulation. This question is answered in \cite{Chen2022}.

First, we introduce a family of optimal control problems.
For any $(t,\eta)\in [0,T)\times H$,   the control  system is
\begin{equation}\label{system2}
\begin{cases}
\ds dX(s)=\big(AX(s)+a(t,X(s),u(s))\big)dt + b(s,X(s),u(s))dW(s), &s \in (t,T],\\
\ns\ds X(t)=\eta,
\end{cases}
\end{equation}
and the cost functional is
\begin{equation}\label{cost2}
    \mathcal{J}(t,\eta;u(\cdot))=\mathbb{E}\Big(\int_t^T
    f(s,X(s),u(s))ds+h(X(T))\Big).
\end{equation}
For any $u(\cdot)\in \mathcal{U}[t,T]$, it follows
immediately from the classical well-posedness of SEEs (e.g., \cite[Theorem 3.14]{Lu2021}) that  the control system
\eqref{system2} has a unique mild solution $X(\cdot) \in C_\dbF([t,T];L^2(\Omega;H))$. Hence, the cost functional \eqref{cost2} is well-defined.

Consider the following optimal control problem:

\textbf{Problem} $\boldsymbol{(S_{t\eta})}$.  For any given $(t,\eta)\in
[0,T]\times H$, find a
$\bar{u}(\cdot)\in \mathcal{U}[t,T]$ such that
\begin{equation}\label{OP2}
    \mathcal{J}(t,\eta;\bar{u}(\cdot))=\inf\limits_{u(\cdot)\in \mathcal{U}[t,T]}\mathcal{J}(t,\eta;u(\cdot)).
\end{equation}

Any $\bar{u}(\cdot)\in \mathcal{U}[t,T]$ satisfying \eqref{OP2} is called an {\it optimal control} (of
\textbf{Problem} $\boldsymbol{(S_{t\eta})}$). The corresponding state
$\overline{X}(\cdot)$ is called an {\it optimal state}, and
$(\overline{X}(\cdot),\bar{u}(\cdot))$ is called an {\it optimal pair}.

The value function is defined as follows:
\begin{equation*}
\begin{cases}
V(t,\eta)=\inf\limits_{u(\cdot)\in \mathcal{U}[t,T]}\mathcal{J}(t,\eta;u(\cdot)),\ \ \forall \ (t,\eta)\in [0,T)\times H,\\
V(T,\eta)=h(\eta),\ \ \forall \ \eta \in H.
\end{cases}
\end{equation*}

It is easy to prove that the value function enjoys the following
properties:
\begin{proposition}\label{prop3.4-1}\cite[Proposition 3.1]{Chen2022} For each $t\in [0,T]$, $\eta$ and
$\eta^{\prime}\in H$, we have
\begin{equation}\label{prop3.4-1-eq1}
|V(t,\eta)|\leq C(1+|\eta|_H)
\end{equation}
and
\begin{equation}\label{prop3.4-1-eq2}
|V(t,\eta)- V(t,\eta^{\prime})| \leq C|\eta-\eta^{\prime}|_H.
\end{equation}
\end{proposition}

\begin{proposition}\label{prop3.6-1}\cite[Proposition 3.2]{Chen2022}
The function $V(\cd,\eta)$ is continuous.
\end{proposition}

We have the following  Dynamic Programming Principle (\cite[Theorem
3.1]{Chen2022}):

\ss

For any $(t,\eta)\in [0,T)\times H$,
\begin{equation}\label{12.15-eq6}
V(t,\eta)
=\inf_{u(\cdot)\in\mathcal{U}[t,T]}\mathbb{E}\Big(\int_t^{\hat
t}f(s, X(s;t,\eta,u),u(s))ds +V\big(\hat t,X(\hat
t;t,\eta,u)\big)\Big),\q\forall\,0\le t\le\hat t \le T.
\end{equation}
Here $X(\cdot;t,\eta,u)$ is the mild solution of \eqref{system2}. By
\eqref{12.15-eq6}, one can derive the HJB equation satisfied by
$V(\cd,\cd)$. We do not present that here since we do not use it in
this paper.

In this paper, we will investigate the
sufficient optimality condition  ---
verification theorem of the \textbf{Problem}
$\boldsymbol{(S_{\eta})}$)---via the  value
function. The verification theorem provides a
way of testing whether a given admissible
control is optimal and enables one to construct
an optimal control via the value function. This
theorem was first studied in the 1960s by
Pontryagin and his group for the LQ problem of
control systems governed by ordinary
differential equations (e.g.,
\cite{Pontryagin1962}). The general cases for
controlled ordinary differential equations were
studied in the sequel by Fleming and Rishel with
smooth value function in \cite{Fleming1975}, and
by Zhou in \cite{Zhou1993} under
viscosity-solution framework. The succedent work
for control systems governed by stochastic
differential equations were studied in
\cite{Li1997,Gozzi2005,Gozzi2010}. When the
value function belongs to $C^{1,2}([0,T]\times
H)$, the verification theorem for
\textbf{Problem} $\boldsymbol{(S_{\eta})}$)
follows similar standard results for the
finite-dimensional case (e.g., \cite[Section
2.5]{Fabbri2017}). However, it is well known
that the value function does not belong to
$C^{1,2}([0,T]\times H)$ in general. This leads
to the study of the verification theorem for
\textbf{Problem} $\boldsymbol{(S_{\eta})}$) for
nonsmooth value function. Along this line, there
are many works when the diffusion term of the
control system is independent of the control
variable (see \cite{Da
Prato2002,Fabbri2017,Federico2009,Federico2018,Fuhrman2006,Fuhrman2002}
and the rich references therein). As far as we
know, there is no published work addressing the
verification theorem for \textbf{Problem}
$\boldsymbol{(S_{\eta})}$ with nonsmooth value
function and control dependent diffusion term.
This does not mean that such problem is not
important. Indeed, the control dependent
diffusion term reflects that the control would
influence the scale of uncertainty, which is
indeed the case in many practical systems. In
this paper, we investigate such problem and
prove the following result.

\begin{theorem}\label{th7.1}
Let Assumptions {\rm ({\bf S1})}--{\rm ({\bf
S2})} hold. Let $V\in C([0,T]\times H)$ be the
value function of \textbf{Problem}
$\boldsymbol{(S_{\eta})}$. Let $\eta\in  H$ be
fixed, and $(\overline X(\cdot),\Bar{u}(\cdot))$
be an admissible pair of \textbf{Problem}
$\boldsymbol{(S_{\eta})}$. Suppose
\begin{equation}\label{13.2-eq01}
A\overline X(\cdot)\in L^2_\dbF(0,T;H),
\end{equation}
and for any $\d>0$, there exists
a constant $C_\d>0$ such that
\begin{equation}\label{13-eq2.6}
|V(t_1,\eta)-V(t_0,\eta)|\leq
C_\d(1+|\eta|^2_H)|t_1-t_0|,\ \ \forall\, t_1,\,
t_0\in [0,T-\d),\ \eta \in H.
\end{equation}

If there exists a triple $(\overline{R},\Bar{p},\overline{P})\in
L_{\mathbb{F}}^2(0,T;\mathbb{R})\times L_{\mathbb{F}}^2(0,T;H)\times
L_{S,\dbF}^2(0,T;\dbS(H))$ such that
\begin{equation}\label{13.2-eq02}
(\overline{R},\Bar{p},\overline{P})\in D_{t+,x}^{1,2,+}V(t,\overline
X(t)), \qquad \mbox{a.e. } (t,\om)\in [0,T]\times \Om
\end{equation}
and
\begin{equation}\label{13.2-eq0}
\mathbb{E}\int_0^T \big(\overline{R}(t)+ \lan  \bar{p}(t),A
\overline X(t)\ran_H +G(t,\overline
X(t),\Bar{u}(t),\Bar{p}(t),\overline{P}(t))\big)dt\leq 0,
\end{equation}
where
\begin{equation}\label{prop3.7-eq2}
G(t,\eta,\rho,p,P)=\frac{1}{2}\lan
Pb(t,\eta,\rho),b(t,\eta,\rho)\ran_{\mathcal{L}_2^0}
+ \lan p, a(t,\eta,\rho)\ran_H-f(t,\eta,\rho),
\end{equation}
$$\hspace{5cm}\forall(t,\eta,\rho,p,P)\in[0,T]\times H\times U\times H\times\dbS(H),$$
then $(\overline X(\cdot),\Bar{u}(\cdot))$ is an optimal pair of
\textbf{Problem} $\boldsymbol{(S_{\eta})}$.
\end{theorem}
\begin{remark}
Theorem \ref{th7.1} is expressed in terms of
superdifferential. It is natural to expect that
a similar result holds for subdifferential.
Unfortunately, the answer is no even for
$H=\dbR$ (see \cite[Example 5.6]{Yong1999} for
example).
\end{remark}

The sufficient condition \eqref{13.2-eq0} can be
replaced by an equivalent condition which looks
much stronger.

\begin{proposition}\label{cor3.2}
Condition \eqref{13.2-eq0} in Theorem
\ref{th7.1} is equivalent to the following:
\begin{eqnarray}\label{13.2-eq07}
& & \overline{R}(t)+ \lan  \bar{p}(t),A
\overline X(t)\ran_H +G(t,\overline
X(t),\Bar{u}(t),\Bar{p}(t),\overline{P}(t))\\
& = &  \overline{R}(t)+ \lan  \bar{p}(t),A
\overline X(t)\ran_H + \inf\limits_{u\in
U}G(t,\overline X(t),
u,\Bar{p}(t),\overline{P}(t)) ,\qquad \mbox{a.e.
} t\in [0,T],\ \mathbb{P}\mbox{-a.s.} \nonumber
\end{eqnarray}
\end{proposition}

Two assumptions, i.e., \eqref{13.2-eq01} and
\eqref{13-eq2.6} are given in Theorem
\ref{th7.1}. The first one is set for the
regularity for the solution of the control
system.  This depends on the SPDE  which governs
the control system. The second one, i.e.,
\eqref{13-eq2.6}, is for the regularity property
of the value function.    This is very subtle
since the value function is not easy to be
handled. In fact, in stochastic case, generally
the value function associated to a control
system is not Lipschitz in $t$ even when all the
coefficients involved are smooth (e.g.
\cite{Buckdahn2010}). Fortunately, our next
result concludes that \eqref{13-eq2.6} holds
under suitable conditions.

\ss

{\bf (S1)$'$} {\it  Suppose that: i) $a(\cdot, \cdot,\cdot): [0,T]\times H\times U \to H$ is
    $\mathcal{B}([0,T])\otimes \mathcal{B}(H)\otimes
    \mathcal{B}(U)/\mathcal{B}(H)$-measurable and
    $b(\cdot,\cdot,\cdot): [0,T]\times H\times U \to
    \mathcal{L}_2^0$ is $\mathcal{B}([0,T])\otimes \mathcal{B}(H)\otimes
    \mathcal{B}(U)/\mathcal{B}(\mathcal{L}_2^0)$-measurable;  ii) for any
    $t\in [0,T],\ \eta\in  H$, the maps $a(t,\eta,\cdot):U\to H $ and
    $b(t,\eta,\cdot): U\to \mathcal{L}_2^0$ are continuous; and iii)
    for any $(t,t_1,t_2,\eta_1,\eta_2,u)\in [0,T]\times [0,T]\times [0,T]\times H \times H\times U$,
\begin{equation*}
\begin{cases}
\big|a(t_1,\eta_1,u)-a(t_2,\eta_2,u)\big|_H
\leq C\big(|t_1-t_2|+|\eta_1-\eta_2|_H\big),\\
\big|b(t_1,\eta_1,u)- b (t_2,\eta_2,u)\big|_{\mathcal{L}_2^0}\leq C\big(|t_1-t_2|+|\eta_1-\eta_2|_H\big),\\
|a(t,0,u)|_H \leq C, \q \ \ \ \
|b(t,0,u)|_{\mathcal{L}_2^0}\leq C.
\end{cases}
\end{equation*}}

\ss

{\bf (S2)$'$} {\it  Suppose that: i) $f(\cdot,\cdot,\cdot):[0,T]\times H\times U\to \mathbb{R}$
    is $\mathcal{B}([0,T])\otimes \mathcal{B}(H)\otimes
    \mathcal{B}(U)/\mathcal{B}(\mathbb{R})$-measurable and
    $h(\cdot):H\to \mathbb{R}$ is
    $\mathcal{B}(H)/\mathcal{B}(\mathbb{R})$-measurable;  ii) For any
    $t\in [0,T],\ \eta\in  H$, the functional $f(t,\eta,\cdot): U\to
    \mathbb{R}$ is continuous; and iii) For any $(t,t_1,t_2,\eta_1,\eta_2,u)\in [0,T]\times [0,T]\times
    [0,T]\times H\times H\times U$,
    \begin{equation*}
        \begin{cases}
            \big|f(t_1,\eta_1,u)-f(t_2,\eta_2,u)\big|\leq C\big(|t_1-t_2|+|\eta_1-\eta_2|_H\big), \\
            \big|h(\eta_1)-h(\eta_2)\big|\leq C|\eta_1-\eta_2|_H\\
            \big|f(t,0,u)\big|\leq C,  \q \ \ |h(0)|\leq C.
        \end{cases}
\end{equation*}}
\begin{theorem}\label{th4.1}
Under Assumptions {\bf (S1)$'$}--{\bf (S2)$'$},
the  value function $V$ is Lipschitz continuous
in $[0,T-\delta]\times H$ for all $\delta
>0$, provided that $A$ generates an analytic semigroup on $H$.
\end{theorem}
\begin{remark}
Generally speaking, the value function is not
Lipschitz continuous in $[0,T]\times H$ even for
$H=\dbR$(e.g.,\cite{Buckdahn2010}).
\end{remark}
Theorem \ref{th4.1} illustrates that the strong
restrictions on $V$ is valid for several
important control systems, in particular, those
governed by stochastic parabolic equations.
What's more, with the additional assumption that
$A$ is analytic, we can drop the assumption
\eqref{13-eq2.6}. Inequality \eqref{13.2-eq4.1}
can be deduced obviously from Theorem
\ref{th4.1} and \eqref{4-eq13}. Therefore, we
have the following result.

\begin{corollary}\label{cor3.1}
Let Assumptions {\rm {\bf (S1)$'$}--{\bf (S2)$'$}} hold. Let $A$ generate
an analytic semigroup, and $V\in C([0,T]\times H)$ be the value
function of \textbf{Problem} $\boldsymbol{(S_{\eta})}$. Let
$(\overline X(\cdot),\Bar{u}(\cdot))$ be an admissible pair of
\textbf{Problem} $\boldsymbol{(S_{\eta})}$. Suppose that
\begin{equation}\label{13.2-eq01.1}
A\overline X(\cdot)\in L^2_\dbF(0,T;H).
\end{equation}
If there exists a triple
$(\overline{R},\Bar{p},\overline{P})\in
L_{\mathbb{F}}^2(0,T;\mathbb{R})\times
L_{\mathbb{F}}^2(0,T;H)\times
L_{S,\mathbb{F}}^2(0,T;\dbS(H))$ such that
\begin{equation}\label{13.2-eq02.1}
\big(\overline{R},\Bar{p},\overline{P}\big)\in
D_{t+,x}^{1,2,+}V(t,\overline X(t)), \qquad
\mbox{a.e. } t\in [0,T],\ \mathbb{P}\mbox{-a.s.}
\end{equation}
and
\begin{equation}\label{13.2-eq0.1}
\mathbb{E}\int_s^T \big(\overline{R}(t)+ \lan
\bar{p}(t),A \overline X(t)\ran_H +G(t,\overline
X(t),\Bar{u}(t),\Bar{p}(t),\overline{P}(t))\big)dt\leq
0,
\end{equation}
then $(\overline X(\cdot),\Bar{u}(\cdot))$ is an
optimal pair of \textbf{Problem}
$\boldsymbol{(S_{\eta})}$.
\end{corollary}

The rest of this paper is divided into three
sections. Section \ref{sec-SVT} is devoted to
the proof of Theorem \ref{th7.1} and Section
\ref{sec-lipschitz} is addressed to  the proof
of Theorem \ref{th4.1}. At last, in Section
\ref{sec-exam}, we provide an illustrative
example fitting for the assumptions in  Theorem
\ref{th7.1}.

\section{Stochastic Verification Theorem}\label{sec-SVT}


In this section, we are going to prove the well-known verification theorem for the infinite dimensional
stochastic control system. The main idea comes from
\cite{Gozzi2005,Yong1999}.

\ss

Let us first recall the concept of regular
conditional probability, which allows us to
regard the conditional expectation as merely
mathematical expectation  taken with respect to
the conditional measure. More details can be
found in \cite[Chapter V, Section
8]{Parthasarathy2005}.

\begin{lemma}\label{prop-2.2}
Let $\mathcal{G}$ be a sub-$\sigma$-algebra of $\mathcal{F}$. Then there exists a map $\mathbf P:  \Omega\times \mathcal{F}\to [0,1]$, called a regular conditional probability given $\mathcal{G}$, such that

\ss

(i)  for each $\omega\in \Omega$, $\mathbf P(\omega,\cdot)$ is a probability measure on $\mathcal{F}$;

\ss

(ii)  for each $A\in \mathcal{F}$, the function $\mathbf P(\cdot,A)$ is $\mathcal{G}$-measurable;

\ss

(iii)  for each $B\in \mathcal{F}$, $\mathbf P(\omega,B)=\mathbb{P}(B|\mathcal{G})(\omega)=\mathbb{E}(1_B|\mathcal{G})(\omega),\ \mathbb{P}${\rm-a.s.}

\ss

We write $\mathbb{P}(\cdot|\mathcal{G})(\omega)$ for $p(\omega,\cdot)$.
\end{lemma}

\ss

The next proposition is taken from
\cite{Yong1999}  with a slight modification.

\begin{proposition}\label{prop-7.1}
Let $v\in C([0,T]\times H)$ and $(t_0,x_0)\in [0,T)\times H$ be given.
Then $(q,p,P)\in D_{t+,x}^{1,2,+}v(t_0,x_0)$ if and only if there
exists a function $\varphi\in C^{1,2}([0,T]\times H)$ such that
\begin{equation}\label{13.1-eq1}
\begin{cases}\ds
\big(\varphi(t_0,x_0),\varphi_t(t_0,x_0),\varphi_x(t_0,x_0),\varphi_{xx}(t_0,x_0)\big)=\big(v(t_0,x_0),q,p,P\big),\\
\ns\ds \varphi(t,x)>v(t,x), \qquad \forall \ (t_0,x_0)\neq (t,x)\in
[t_0,T]\times H.
\end{cases}
\end{equation}
\end{proposition}

{\it Proof}. The ``if" part follows directly from the definition of
$D_{t+,x}^{1,2,+}v(t_0,x_0)$.

The ``only if" part. Suppose $(q,p,P)\in
D_{t+,x}^{1,2,+}v(t_0,x_0)$. Define a functional
on $[t_0,T]\times H$ as follows: if
$(t_0,x_0)\neq (t,x) \in [t_0,T]\times H$, then
\begin{equation*}
\Phi(t,x)=
\Big[\frac{v(t,x)-v(t_0,x_0)-q(t-t_0)-\lan
p,x-x_0\ran_H-\frac{1}{2}\lan
P(x-x_0),x-x_0\ran_H}{t-t_0+|x-x_0|_H^2}\Big]^+;
\end{equation*}
and if $(t_0,x_0)\neq (t,x)$, then
$\Phi(t,x)=0$.

Let
\begin{equation*}
\k(r)= \begin{cases}\ds\sup \big\{\Phi(t,x)\big|
(t,x)\in(t_0,T]\times H,\ t-t_0+|x-x_0|_H^2\leq
r\big\}, &\mbox{ if }
r>0,\\
\ns\ds 0, &\mbox{ if }r\leq 0.
\end{cases}
\end{equation*}
Then $\k:\mathbb{R}\to [0,+\infty)$ is a
continuous and nondecreasing function with
$\k(0)=0$. Further, we have
\begin{eqnarray*}
& & v(t,x)-\Big[v(t_0,x_0)+q(t-t_0)+\lan p,x
-x_0\ran_H+\frac{1}{2} \lan
P(x-x_0),x-x_0\ran_H\Big]
\\
& \leq\3n &
\big(t-t_0+|x-x_0|_H^2\big)\k\big(t-t_0+|x-x_0|_H^2\big),\qquad\qq
\forall \ (t,x)\in [t_0,T]\times H.
\end{eqnarray*}
Let
$$
\Psi(\rho)=\frac{2}{\rho}\int_0^{2\rho}
\int_0^r\k(\theta)d\theta dr,\qquad \rho>0.
$$
Then we have
$$
\Psi_{\rho}(\rho)=-\frac{2}{\rho
^2}\int_0^{2\rho} \int_0^r\k(\theta)d\theta dr +
\frac{4}{\rho}\int_0^{2\rho}\k(\theta)d\theta
$$
and
$$\Psi_{\rho\rho}(\rho)
=\frac{4}{\rho ^3}
\int_0^{2\rho}\int_0^r\k(\theta)d\theta dr-
\frac{8}{\rho
^2}\int_0^{2\rho}\k(\theta)d\theta+
\frac{8}{\rho}\k(2\rho).$$
Consequently,
$$|\Psi(\rho)|\leq 4\rho\k(2\rho),
\qquad |\Psi_{\rho}(\rho)|\leq
12\k(2\rho),\qquad |\Psi_{\rho\rho}(\rho)|\leq
\frac{32\k(2\rho)}{\rho}. $$
Now we define
\begin{equation}\label{12.14-eq3}
\psi(t,x)=\begin{cases}\ds \Psi(\rho(t,x))+\rho(t,x)^2, &\mbox{ if }
(t_0,x_0)\neq (t,x)\in
[t_0,T]\times H,\\
\ns\ds 0, &\mbox{ if } (t,x)=(t_0,x_0),
\end{cases}
\end{equation}
where $\rho(t,x)=t-t_0+|x-x_0|_H^2$. Set
\begin{eqnarray}\label{13.1-eq2}
\varphi(t,x) \3n & =\3n & v(t_0,x_0)+q(t-t_0)+\lan p,x-x_0\ran_H\nonumber\\
& & +\frac{1}{2}\lan P(x-x_0),x-x_0\ran_H + \psi(t,x),\qquad \forall
\ (t,x)\in [0,T]\times H.
\end{eqnarray}
We claim that $\varphi\in C^{1,2}([0,T]\times H)$ satisfies
\eqref{13.1-eq1}. First, for any $(t,x)\in [t_0,T]\times H$ with
$(t,x)\neq (t_0,x_0)$, we have
\begin{eqnarray*}
\psi(t,x)\3n& >\3n & \frac{2}{\rho(t,x)}
\int_{\rho(t,x)}^{2\rho(t,x)}\int_0^r\k(\theta)d\theta
dr
\\
\3n& \ge\3n & \frac{2}{\rho(t,x)}\k(\rho(t,x))
\int_{\rho(t,x)}^{2\rho(t,x)}\big(r-\rho(t,x)\big)dr
\\
\3n& = \3n& \rho(t,x)\k(\rho(t,x)).
\end{eqnarray*}
Next,  for any $(t,x)\in [t_0,T]\times H$, it follows from
\eqref{12.14-eq3} that
$$\psi_t(t,x)=\Psi_{\rho}(\rho(t,x))+2\rho(t,x),$$
$$\psi_{x}(t,x)=2\Psi_{\rho}(\rho(t,x))(x-x_0)+4\rho(t,x)(x-x_0),$$
and
$$\psi_{xx}(t,x)=4\Psi_{\rho\rho}(\rho(t,x))\otimes(x-x_0)+2\Psi_{\rho}(\rho(t,x))I+4\rho(t,x)I+8(x-x_0)\otimes(x-x_0).$$
Thus,  noting $|x-x_0|\leq \rho(t,x)$, we obtain
\begin{equation*}
\begin{cases}\ds
|\psi(t,x)|\leq 4\rho(t,x)\k(2\rho(t,x))+\rho(t,x)^2,\\
\ns\ds |\psi_t(t,x)\leq 12\k(\rho(t,x))+2\rho(t,x),\\
\ns\ds |\psi_x(t,x)\leq 24|x-x_0|\k(2\rho(t,x))+4\rho(t,x)|x-x_0|,\\
\ns\ds |\psi_{xx}(t,x)|\leq
\frac{128|x-x_0|^2}{\rho(t,x)}\k(2\rho(t,x))+24\k(2\rho(t,x))+12\rho(t,x)\leq
152\k(2\rho(t,x))+12\rho(t,x).
\end{cases}
\end{equation*}
Hence, $\psi\in C^{1,2}([0,T]\times H)$ and
$$
\psi(t_0,x_0)=0, \qquad \psi_t(t_0,x_0)=0,\qquad \psi_x(t_0,x_0)=0,
\qquad \psi_{xx}(t_0,x_0)=0.
$$
This proves our claim.
\endpf

To continue, we need two known results.   The
first one is taken from from \cite{Yong1999} but
with a slight modification following the
discussion in \cite{Gozzi2010} and \cite[Remark
3.4, Section 3]{Federico2009}.
\begin{lemma}\label{lm7.1}
Let $g\in C[0,T]$. Extend $g$ to
$(-\infty,+\infty)$ with $g(t)=g(T)$, for $t>T$,
and $g(t)=g(0)$ for $t<0$.  Suppose that for
each $\d\in(0,T)$, there is a $\rho_\d\in
L^1(0,T-\d)$ such that for some $\e_0>0$,
$$ \frac{g(t+\e)-g(t)}{\e}\leq \rho(t), \qquad
\forall \e\leq\e_0,\  \ae \ t\in [0,T-\d).$$
Then
$$
g(\beta)-g(\alpha)\leq \int_{\alpha}^{\beta}
\uplim_{\e\to 0^+}\frac{g(r+\e)-g(r)}{\e}dr,
\qquad \forall\, 0\leq \alpha\leq \beta\leq
T-\d.
$$
\end{lemma}
\begin{lemma}\label{lm7.2}\cite[Theorem 9, Chapter
2]{Diestel1977} Let $Z$ be a Banach space, $[a,b]\subset \mathbb{R}$
and $z:[a,b]\to Z$ be a Bochner integrable function. Then
\begin{equation*}
\frac{1}{\e}\int_t^{t+\e}|z(r)-z(t)|_Zdr\to 0,\
\mbox{\rm as} \ \e\to 0^+,\ \ \ae \ t\in [a,b].
\end{equation*}
\end{lemma}

According to \cite[page 92]{Doob1994}, if $(\Om,
\cF, \dbP)$ is a complete probability space, $Z$
is a separable Hilbert space, and $\cF$ is
countably generated apart from null sets,  then
$L^1(\Omega; Z)$ is separable. Hence,
$L^1(\Omega;H)$ is separable Banach space. Then,
following \cite[Lemma 3.6, Section
3.2]{Gozzi2005}, we have an analogous result.
\begin{lemma}\label{lm7.3}
Let $z\in L^1_{\mathbb{F}}(0,T;H)$. Then it is
Bochner integrable if it is regarded as a map
from $[0,T]$ to $L^1(\Omega;H)$.
\end{lemma}

Now we are in a position to prove Theorem \ref{th7.1}.

\ss

{\it Proof of Theorem \ref{th7.1}}. We divide
the proofs into several steps.

{\bf Step 1}. By Assumption \textbf{(S1)} and
\eqref{13.2-eq0},  $z(\cdot)=\bar{a}(\cdot),\
\bar{b}(\cdot), \ A\overline X(\cdot)$ can be
regarded as  Bochner integrable functions from
$[0,T]$ to $L^1(\Omega,H)$ (\cite[Lemma 3.6,
Section 3.2]{Gozzi2005}). Noting that by Lemma
\ref{lm7.2}, for $\e$ small enough,
\begin{equation}\label{8.19-eq1}
\frac{1}{\e}\int_t^{t+\e} |z(r)-z(t)|_{L^1(\Omega,H)}dr< \infty,
\end{equation}
and
\begin{equation}\label{8.19-eq2}
\frac{1}{\e}\int_t^{t+\e} |z(r,\cdot)-z(t,\cdot)|_Hdr\leq C,\q \dbP{-a.s.}
\end{equation}
By \eqref{8.19-eq1}, \eqref{8.19-eq2}, Fubini's Theorem and the
Dominated Convergence Theorem, we have
\begin{equation}\label{13.2-eq03}
\lim_{\e\to
0^+}\mathbb{E}\frac{1}{\e}\int_t^{t+\e}|z(r)-z(t)|_Hdr=
0,\q  \ \mbox{a.e.} \ t\in [0,T].
\end{equation}
Fix $t_0\in [0,T]$ such that \eqref{13.2-eq02} holds at $t_0$, and
\eqref{13.2-eq03} holds at $t_0$ for $z(\cdot)=A\overline X(\cdot)$,
$\bar{a}(\cdot)$ and $\bar{b}(\cdot)$.

Fix $\omega_0\in \Omega$ such that the regular conditional
probability $\mathbb{P}\big(\cdot|\mathcal{F}_{t_0}\big)(\omega_0)$
is well-defined. In the probability space $\big(\Omega, \mathcal{F},
\mathbb{P}\big(\cdot|\mathcal{F}_{t_0}\big)(\omega_0)\big)$, the
random variables
$$\overline X(t_0),\ \overline R(t_0),\ \overline p(t_0),\ \overline R(t_0)$$
are almost surely equal to
$$\overline X(t_0,\omega_0),\ \overline R(t_0,\omega_0),\ \overline p(t_0,\omega_0),\ \overline R(t_0,\omega_0),$$
respectively. Let $\eta_0=\overline X(t_0,\omega_0)$. Denote by
$\mathbb{E}_{\omega_0}$ the expectation with respect to the
probability measure
$\mathbb{P}\big(\cdot|\mathcal{F}_{t_0}\big)(\omega_0)$.

By Proposition \ref{prop-7.1}, there exists  a
function $\phi \in C^{1,2}([0,T]\times H)$ such
that
\begin{equation}\label{13.2-eq1}
\phi(t,\eta)> V(t,\eta), \qquad \mbox{for every } (t,\eta)\in
(0,T)\times H, \ (t,\eta)\neq (t_0,\eta_0),
\end{equation}
\begin{equation}\label{13.2-eq2}
\big(\phi(t_0,\eta_0),\phi_t(t_0,\eta_0),\phi_x(t_0,\eta_0),\phi_{xx}(t_0,
\eta_0)\big)=\big(V(t_0,\eta_0),
\overline{R}(t_0,\omega_0),\Bar{p}(t_0,\omega_0),\overline{P}(t_0,\omega_0)\big),
\end{equation}
and that
\begin{equation}\label{12.14-eq4}
\phi,\ \phi_t,\ \phi_x,\ \phi_{xx} \mbox{ are polynomially bounded}.
\end{equation}
Then $\phi$ is a fixed deterministic function if
$(t_0,\omega_0)$ is fixed.

Applying It\^o's formula to $\phi$, then for any
$\e>0$, it follows that
\begin{eqnarray}\label{12.14-eq14}
& &  \phi\big(t_0+\e,\overline X(t_0+\e)\big)-
\phi\big(t_0,\overline X(t_0)\big)\nonumber\\
\3n& =\3n &
\int_{t_0}^{t_0+\e}\Big(\phi_t\big(r,\overline
X(r)\big)\! + \lan \phi_x\big(r,\overline
X(r)\big),\Bar{a}(r)\ran_H\!+\!\lan
\phi_x\big(r,\overline X(r)\big),A\overline
X(r)\ran_H
\\
& & \! + \frac{1}{2}\lan
\phi_{xx}\big(r,\overline
X(r)\big)\Bar{b}(r),\Bar{b}(r)\ran_{\mathcal{L}_2^0}\Big)dr+\int_{t_0}^{t_0+\e}\lan
\phi_x\big(r,\overline
X(r)\big),\Bar{b}(r)dW(r)\ran_H.\nonumber
\end{eqnarray}
Let $\{\e_n\}_{n=1}^\infty\subset (0,+\infty)$
be such that $\lim\limits_{n\to+\infty}\e_n=0$.
From \eqref{12.14-eq14}, we have
\begin{eqnarray}\label{13.2-eq3}
& &
\mathbb{E}_{\omega_0}\frac{1}{\e_n}\big(V\big(t_0+\e_n,
\overline X(t_0+\e_n)\big)- V\big(t_0,\overline
X(t_0)\big)\big)\nonumber
\\
\3n& \leq\3n &
\mathbb{E}_{\omega_0}\frac{1}{\e_n}
\int_{t_0}^{t_0+\e_n}
\Big(\phi_t\big(r,\overline X(r)\big)+
\lan\phi_x\big(r,\overline X(r)\big),A\overline
X(r)\ran_H+\lan
\phi_x\big(r,\overline X(r)\big),\Bar{a}(r)\ran_H\\
& & \hspace{2.35cm}+\frac{1}{2}\lan
\phi_{xx}\big(r,\overline
X(r)\big)\Bar{b}(r),\Bar{b}(r)\ran_{\mathcal{L}_2^0}\Big)dr.
\nonumber
\end{eqnarray}
\ss

{\bf Step 2}. In this step,  we  treat  the
right hand side of \eqref{13.2-eq3} term by
term.

\ss

First, thanks to the continuity of $\overline X(\cdot)$ and
$\phi_t(\cdot)$, we get
\begin{equation}\label{12.14-eq15}
\lim_{\e_n\to
0^+}\frac{1}{\e_n}\int_{t_0}^{t_0+\e_n}\phi_t\big(r,\overline
X(r)\big)dr = \phi_t\big(t_0,\overline
X(t_0)\big),\qquad
\mathbb{P}(\cdot|\mathcal{F}_{t_0})
(\omega_0)\mbox{-a.s.}
\end{equation}
By the polynomial growth of $\phi_t$, due to the
Dominated Convergence Theorem, it holds that
\begin{equation}\label{12.14-eq15.1}
\lim_{\e_n\to
0^+}\mathbb{E}_{\omega_0}\Big|\frac{1}{\e_n}
\int_{t_0}^{t_0+\e_n}\phi_t\big(r,\overline
X(r)\big)dr- \phi_t\big(t_0,\overline
X(t_0)\big)\Big|_H=0.
\end{equation}

As for the second term of \eqref{13.2-eq3}, we
have
\begin{eqnarray}\label{13.2-eq4}
& & \mathbb{E}_{\omega_0}\frac{1}{\e_n}
\int_{t_0}^{t_0+\e_n}\lan \phi_x\big(r,
\overline X(r)\big),\Bar{a}(r)\ran_H dr-\lan
\phi_x\big(t_0,\overline
X(t_0)\big),\bar{a}(t_0)\ran_H \nonumber
\\
\3n& = \3n&
\mathbb{E}_{\omega_0}\frac{1}{\e_n}\int_{t_0}^{t_0+\e_n}\lan
\phi_x\big(r,\overline
X(r)\big)-\phi_x\big(t_0,\overline
X(t_0)\big),\Bar{a}(r)\ran_H dr
\\
& &
+\mathbb{E}_{\omega_0}\frac{1}{\e_n}\int_{t_0}^{t_0+\e_n}\lan
\phi_x\big(t_0,\overline
X(t_0)\big),\Bar{a}(r)-\bar{a}(t_0)\ran_H dr.
\nonumber
\end{eqnarray}
Clearly,
\begin{eqnarray}\label{12.14-eq10}
& &
\Big|\mathbb{E}_{\omega_0}\frac{1}{\e_n}\int_{t_0}^{t_0+\e_n}\lan
\phi_x\big(r,\overline
X(r)\big)-\phi_x\big(t_0,\overline
X(t_0)\big),\Bar{a}(r)\ran_H dr\Big|\nonumber
\\
\3n& \leq \3n& \mathbb{E}_{\omega_0}
\frac{1}{\e_n}\int_{t_0}^{t_0+\e_n}\big|
\phi_x\big(r,\overline
X(r)\big)-\phi_x\big(t_0,\overline
X(t_0)\big)\big|_H\big|\Bar{a}(r)\big|_H dr
\\
\3n& \leq \3n& \mathbb{E}_{\omega_0}
\Big(\frac{1}{\e_n}\int_{t_0}^{t_0+\e_n}\big|
\phi_x\big(r,\overline
X(r)\big)-\phi_x\big(t_0,\overline
X(t_0)\big)\big|_H^2dr\Big)^{\frac{1}{2}}
\Big(\frac{1}{\e_n}\int_{t_0}^{t_0+\e_n}\big|\Bar{a}(r)\big|_H^2
dr\Big)^{\frac{1}{2}}\nonumber
\\
\3n& \leq \3n&
\Big(\mathbb{E}_{\omega_0}\frac{1}{\e_n}\int_{t_0}^{t_0+\e_n}\big|
\phi_x\big(r,\overline
X(r)\big)-\phi_x\big(t_0,\overline
X(t_0)\big)\big|_H^2dr\Big)^{\frac{1}{2}}
\Big(\mathbb{E}_{\omega_0}\frac{1}{\e_n}
\int_{t_0}^{t_0+\e_n}\big|\Bar{a}(r)\big|_H^2
dr\Big)^{\frac{1}{2}}.\nonumber
\end{eqnarray}
By Assumption {\bf (S1)}, it follows that
\begin{eqnarray}\label{12.14-eq5}
\mathbb{E}_{\omega_0}\frac{1}{\e_n}
\int_{t_0}^{t_0+\e_n}\big|\Bar{a}(r)\big|_H^2 dr  \leq
C\frac{1}{\e_n}\int_{t_0}^{t_0+\e_n}
\big(1+\mathbb{E}_{\omega_0}\big|\overline X(r)\big|^2_H\big)dr \leq
C\big(1+|\eta|_H^2\big).
\end{eqnarray}
Arguing as in Step 1 for $\phi_t$, we get that
\begin{eqnarray}\label{12.14-eq6}
\lim_{n\to
+\infty}\mathbb{E}_{\omega_0}\frac{1}{\e_n}\int_{t_0}^{t_0+\e_n}
\big|\phi_x\big(r,\overline
X(r)\big)-\phi_x\big(t_0,\overline
X(t_0)\big)\big|_H dr=0.
\end{eqnarray}
By \eqref{12.14-eq10}--\eqref{12.14-eq5}, we see the first term on
the right hand side of \eqref{13.2-eq4} goes to zero as $n\to
+\infty$.

Now we handle the second term of
\eqref{13.2-eq4}. Clearly,
\begin{eqnarray}\label{12.14-eq11}
& & \Big|\mathbb{E}_{\omega_0}\frac{1}{\e_n}
\int_{t_0}^{t_0+\e_n} \lan
\phi_x\big(t_0,\overline
X(t_0)\big),\Bar{a}(r)-\bar{a}(t_0)\ran_H
dr\Big|\nonumber
\\
& & \leq \big|\phi_x\big(t_0,\overline
X(t_0)\big)\big|_H\mathbb{E}_{\omega_0}
\frac{1}{\e_n}\int_{t_0}^{t_0+\e_n}\big|\Bar{a}(r)-\bar{a}(t_0)\big|_Hdr.
\end{eqnarray}
By the choice of $t_0$, we have
\begin{eqnarray*}\label{12.14-eq12}
0\3n& = \3n&\lim_{n\to
+\infty}\mathbb{E}\frac{1}{\e_n}\int_{t_0}^{t_0+\e_n}\big|\bar{a}(r)-\bar{a}(t_0)\big|_H dr \nonumber\\
\3n& =\3n& \lim_{n\to
+\infty}\mathbb{E}\Big[\mathbb{E}\Big(\frac{1}{\e_n}\int_{t_0}^{t_0+\e_n}
\big|\bar{a}(r)-\bar{a}(t_0)\big|_Hdr\Big|\mathcal{F}_{t_0}^s\Big)\Big]
\\
\3n& =\3n& \lim_{n\to
+\infty}\mathbb{E}\Big(\mathbb{E}_{\omega_0}\frac{1}{\e_n}\int_{t_0}^{t_0+\e_n}\big|\bar{a}(r)-\bar{a}(t_0)\big|_Hdr\Big).
\nonumber
\end{eqnarray*}
This implies that
\begin{eqnarray*}
\lim_{n\to
+\infty}\mathbb{E}_{\omega_0}\frac{1}{\e_n}
\int_{t_0}^{t_0+\e_n}\big|\bar{a}(r)-\bar{a}(t_0)\big|_Hdr=0
\mbox{ in $L^1_\cF(\Omega;\mathbb{R})$}.
\end{eqnarray*}
Hence, there is a subsequence  $\{\e_n^{(1)}\}_{l=1}^\infty$ of
$\{\e_n\}_{n=1}^\infty$ such that for $\mathbb{P}\mbox{-a.s.}\
\omega_0$,
$$
\lim\limits_{n\to+\infty}
\mathbb{E}_{\omega_0}\frac{1}{\e_n^{(1)}}
\int_{t_0}^{t_0+\e_n^{(1)}}\big|\bar{a}(r)-\bar{a}(t_0)\big|_Hdr=
0.
$$
This, together with \eqref{12.14-eq11}, implies that
\begin{equation}\label{12.14-eq12.1}
\lim\limits_{n\to+\infty}\Big|
\mathbb{E}_{\omega_0}\frac{1}{\e_n^{(1)}}
\int_{t_0}^{t_0+\e_n^{(1)}} \lan
\phi_x\big(t_0,\overline
X(t_0)\big),\Bar{a}(r)-\bar{a}(t_0)\ran_H
dr\Big|=0.
\end{equation}

\ss

Next, we treat the third term of \eqref{13.2-eq3}. Obviously,
\begin{eqnarray}\label{13.2-eq3-1}
& &
\mathbb{E}_{\omega_0}\frac{1}{\e_n}\int_{t_0}^{t_0+\e_n}\lan
\phi_x\big(r, \overline X(r)\big), A\overline
X(r)\ran_Hdr-\lan
\phi_x\big(t_0,\overline X(t_0)\big),A\overline X(t_0)\ran_H \nonumber\\
\3n& = \3n& \mathbb{E}_{\omega_0}\frac{1}{\e_n}
\int_{t_0}^{t_0+\e_n}\lan \phi_x\big(r,\overline
X(r)\big)-\phi_x\big(t_0,\overline
X(t_0)\big),A\overline X(t_0)\ran_H dr\\
& & +\mathbb{E}_{\omega_0}\frac{1}{\e_n}
\int_{t_0}^{t_0+\e_n}\lan \phi_x\big(r,\overline
X(r)\big),A\overline X(r)-A\overline
X(t_0)\ran_H dr.\nonumber
\end{eqnarray}
The first term of \eqref{13.2-eq3-1} reads
\begin{eqnarray*}
& &
\Big|\mathbb{E}_{\omega_0}\frac{1}{\e_n}\int_{t_0}^{t_0+\e_n}\lan
\phi_x\big(r,\overline
X(r)\big)-\phi_x\big(t_0,\overline
X(t_0)\big),A\overline X(t_0)\ran_H dr\Big|\\
\3n& \leq \3n& \mathbb{E}_{\omega_0}
\frac{1}{\e_n}\int_{t_0}^{t_0+\e_n}\big|
\phi_x\big(r,\overline
X(r)\big)-\phi_x\big(t_0,\overline X(t_0)\big)\big|_H\big|A\overline X(t_0)\big|_H dr\\
\3n& \leq \3n& \big|A\overline
X(t_0)\big|_H\mathbb{E}_{\omega_0}
\frac{1}{\e_n}\int_{t_0}^{t_0+\e_n}\big|
\phi_x\big(r,\overline
X(r)\big)-\phi_x\big(t_0,\overline
X(t_0)\big)\big|_Hdr. \nonumber
\end{eqnarray*}
This, together with \eqref{12.14-eq6}, implies
that the first term in \eqref{13.2-eq3-1} tends
to $0$ as $n\to +\infty$. Now we handle the
second term of \eqref{13.2-eq3-1}. Since
$\phi_x(\cdot)$ and  $\overline X(\cdot)$ are
continuous, we see that
\begin{eqnarray}\label{12.14-eq13}
&& \Big|
\mathbb{E}_{\omega_0}\frac{1}{\e_n}\int_{t_0}^{t_0+\e_n}\lan
\phi_x\big(r,\overline X(r)\big),A\overline X(r)-A\overline X(t_0)\ran_H dr\Big|\nonumber\\
& & \leq \Big(\mathbb{E}_{\omega_0}
\frac{1}{\e_n}\int_{t_0}^{t_0+\e_n}\big|\phi_x\big(r,\overline
X(r)\big)\big|_H^2\Big)^{\frac{1}{2}}
\Big(\mathbb{E}_{\omega_0}
\frac{1}{\e_n}\int_{t_0}^{t_0+\e_n}\big|A\overline X(r)-A\overline X(t_0)\big|_H^2dr\Big)^{\frac{1}{2}} \\
& & =\big|\phi_x\big(t_0,\overline
X(t_0)\big)\big|_H^2 \Big(\mathbb{E}_{\omega_0}
\frac{1}{\e_n}\int_{t_0}^{t_0+\e_n}\big|A\overline
X(r)-A\overline
X(t_0)\big|_H^2dr\Big)^{\frac{1}{2}}.\nonumber
\end{eqnarray}
By the choice of $t_0$, following the same
procedure  for deducing \eqref{12.14-eq12.1}, we
get from \eqref{12.14-eq13} that, there exists a
subsequence $\{\e_n^{(2)}\}_{j=1}^\infty$ of
$\{\e_n^{(1)}\}_{n=1}^\infty$  such that
\begin{equation}\label{12.14-eq12-1}
\lim\limits_{n\to\infty}\Big|\mathbb{E}_{\omega_0}
\frac{1}{\e_n^{(2)}}\int_{t_0}^{t_0+\e_n^{(2)}}
\lan \phi_x\big(r,\overline
X(r)\big),\Bar{a}(r)-\bar{a}(t_0)\ran_H
dr\Big|=0.
\end{equation}

\ss

At last, we deal with the forth term of
\eqref{13.2-eq3}. We have
\begin{eqnarray}\label{13.2-eq5-1}
& &
\frac{1}{2}\mathbb{E}_{\omega_0}\Big[\frac{1}{\e_n}
\int_{t_0}^{t_0+\e_n} \big(\lan
\phi_{xx}\big(r,\overline
X(r)\big)\Bar{b}(r),\Bar{b}(r)\ran_{\mathcal{L}_2^0}-
\lan \phi_{xx}\big(t_0,\overline
X(t_0)\big)\Bar{b}(t_0),\Bar{b}(t_0)\ran_{\mathcal{L}_2^0}\big)
dr\Big]\nonumber
\\
\3n& = \3n& \frac{1}{2}\mathbb{E}_{\omega_0}
\Big[\frac{1}{\e_n} \int_{t_0}^{t_0+\e_n}\lan
\big(\phi_{xx}(r,\overline
X(r))-\phi_{xx}(t_0,\overline
X(t_0))\big)\Bar{b}(r),\Bar{b}(r)\ran_{\mathcal{L}_2^0}dr\Big]
\\
& & + \frac{1}{2}\mathbb{E}_{\omega_0}
\Big[\frac{1}{\e_n}\int_{t_0}^{t_0+\e_n}\big(\lan
\phi_{xx}\big(t_0,\overline
X(t_0)\big)\Bar{b}(r),\Bar{b}(r)\ran_{\mathcal{L}_2^0}-\lan
\phi_{xx}\big(t_0,\overline
X(t_0)\big)\Bar{b}(t_0),\Bar{b}(r)\ran_{\mathcal{L}_2^0}\big)dr\Big]\nonumber
\\
& &
+\frac{1}{2}\mathbb{E}_{\omega_0}\Big[\frac{1}{\e_n}\int_{t_0}^{t_0+\e_n}\big(\lan
\phi_{xx}\big(t_0,\overline
X(t_0)\big)\Bar{b}(t_0),\Bar{b}(r)\ran_{\mathcal{L}_2^0}-\lan
\phi_{xx}\big(t_0,\overline
X(t_0)\big)\Bar{b}(t_0),\Bar{b}(t_0)\ran_{\mathcal{L}_2^0}\big)dr\Big].\nonumber
\end{eqnarray}
Now employing the same arguments used to  show
the right-hand side of \eqref{13.2-eq4}
approaching zero, we reach that the right-hand
side of \eqref{13.2-eq5-1} vanishes if we
replace $\{\e_n\}_{n=1}^\infty$ by a subsequence
$\{\e_n^{(2)}\}_{n=1}^\infty$ of
$\{\e_{j}\}_{j=1}^\infty$ and let  $n\to\infty$.

In summary,  for any sequence
$\{\e_{n}\}_{n=1}^\infty\subset (0,+\infty)$
with $\lim\limits_{n\to\infty}\e_n=0$, there
exists a subsequence
$\{\e_n^{(2)}\}_{n=1}^\infty$ of
$\{\e_{n}\}_{n=1}^\infty$, such that
\begin{eqnarray*}
& & \lim_{n\to +\infty}\mathbb{E}_{\omega_0}
\Big[\frac{1}{\e_n^{(2)}}\int_{t_0}^{t_0+\e_n^{(2)}}
\Big(\phi_t\big(r,\overline X(r)\big)+ \lan
\phi_x(r,\overline X(r)),A\overline
X(r)\ran_H+\lan \phi_x\big(r,\overline
X(r)\big),\Bar{a}(r)\ran_H
\\
& & \hspace{4.12cm}+\frac{1}{2}\lan
\phi_{xx}\big(r,\overline X(r)\big)
\Bar{b}(r),\Bar{b}(r)\ran_{\mathcal{L}_2^0}\Big)dr\Big]
\\
&&   = \phi_t(t_0,\overline X(t_0)) + \lan
\phi_x\big(t_0,
\overline{X}(t_0)\big),A\overline X(t_0)\ran_H +
\lan \phi_x(t_0,\overline
X(t_0)),\Bar{a}(t_0)\ran_H \\
&&  \q + \frac{1}{2}\lan
\phi_{xx}\big(t_0,\overline
X(t_0)\big)\Bar{b}(t_0),\Bar{b}(t_0)\ran_{\mathcal{L}_2^0}.
\end{eqnarray*}

\ss

{\bf Step 3}. In this step, we are to prove the
following claim:

\ss

{\bf Claim 1}: {\it For any $\d\in (0,T)$, there
exists  $\rho_\d(\cdot)\in L^1(0,T-\d)$ such
that for almost every $t_0\in [0,T-\d)$ chosen
as in Step 1 and $\e>0$ with $t_0+\e\leq T-\d$,
$\eta\in H$ and $\big(\overline
X(\cdot),\bar{u}(\cdot)\big)$ being the
admissible pair, it holds that
\begin{equation}\label{13.2-eq4.1}
\frac{1}{\e}\mathbb{E}\left(V\big(t_0+\e,\overline
X(t_0+\e)\big)-V\big(t_0,\overline
X(t_0)\big)\right)\leq \rho(t_0).
\end{equation}
}

\ss

Indeed, by Assumption  \eqref{13-eq2.6} and the
definition of $D_{t,x}^{1,2,+}v(\cdot,\cdot)$,
we have
\begin{eqnarray}\label{12-eq2.6}
&&  V\big(t_0+\e,\overline X(t_0+\e)\big)-
V\big(t_0,\overline X(t_0)\big) \\
&&  \leq  C_\d\e\big(1+|\overline
X(t_0+\e)|_H^2\big)+\lan p(t_0), \overline
X(t_0+\e)-\overline X(t_0)\ran_H
+C_0\big|\overline X(t_0+\e)-\overline
X(t_0)\big|_H^2 \nonumber
\end{eqnarray}
and
$$|\bar{p}(t)|_H\leq C\big(1+|\overline X(t)|_H^2\big),\qq \forall t\in [0,T].$$
Here and in what follows, we use $C_\d$ to
denote a constant depending on $\d$, which may
vary from line to line.

Now, we begin to estimate the right hand side of
\eqref{12-eq2.6} term by term.   Noting the
choice of $t_0$, we have
\begin{eqnarray}\label{12-eq2.12}
& & \mathbb{E}\lan \bar{p}(t_0), \overline X(t_0+\e)-\overline X(t_0)\ran_H\nonumber\\
\3n& = \3n& \mathbb{E} \Big\langle \bar{p}(t_0), (S(\e)-I)\overline X(t_0)+\int_{t_0}^{t_0+\e} S(r-t)a\big(r,\overline X(r),\bar{u}(r)\big)dr\nonumber\\
& & \hspace{1.3cm}+\int_{t_0}^{t_0+\e} S(r-t)b\big(r,\overline X(r),\bar{u}(r)\big)dW(r)\Big\rangle_H\\
\3n& \leq \3n& \e \mathbb{E} \lan \bar{p}(t_0),
A\overline X(t_0)\ran_H+ C\big[\mathbb{E}
\big(1+|\overline
X(t_0)|_H^2\big)^2\big]^{1/2}\Big(\mathbb{E}\Big|
\int_{t_0}^{t_0+\e}a\big(r,\overline
X(r),\bar{u}(r)\big)dr\Big|_H^2\Big)^{1/2}
\nonumber\\
\3n& \leq \3n&
C\e\Big[\mathbb{E}\big(1+\big|\overline
X(t_0)\big|_H^2\big)\Big]^{1/2}\Big(\mathbb{E}\big|A\overline
X(t_0)\big|^2_H\Big)^{1/2}+C\e\Big[\mathbb{E}\big(1+\big|\overline
X(t_0)\big|_H^2\big)\Big]^{1/2}.\nonumber
\end{eqnarray}
The third term in \eqref{12-eq2.6} reads
\begin{eqnarray}\label{12-eq2.13}
\mathbb{E}\big|\overline X(t+\e)- \overline
X(t)\big|_H^2 \3n& \leq \3n&
C\Big[\mathbb{E}\big|\big(S(\e)-I\big)\overline X(t)\big|_H^2 + \mathbb{E}\Big| \int_t^{t+\e}a\big(r,\overline X(r),\bar{u}(r)\big)dr\Big|_H^2 \nonumber\\
& & \q +
\mathbb{E}\Big|\int_t^{t+\e}b\big(r,\overline
X(r),\bar{u}(r)\big)dW(r)\Big|_H^2\Big]
\\
\3n& \leq \3n& C\left[\mathbb{E}|A\overline
X(t)|_H^2 \e^2 + \mathbb{E}\big(1+|\overline
X(t)|_H^2\big)(\e^2+\e)\right].\nonumber
\end{eqnarray}
Thus, by taking
$$
\rho(t_0)=C_\d\Big[\mathbb{E}\big(1+ |\overline
X(t_0)|_H^2\big)\Big]^{1/2}
\Big[\Big(\mathbb{E}|A\overline
X(t_0)|^2_H\Big)^{1/2}+1\Big] \in L^1(0,T-\d),$$
we complete the proof of {\bf Claim 1}.

\ss

{\bf Step 4}. Applying {\bf Claim 1} shown in
Step 3 on
 $\big(\Omega, \mathcal{F},
\mathbb{P}\big(\cdot|\mathcal{F}_{t_0}\big)(\omega_0)\big)$,
then by Lemma \ref{lm7.1} and \eqref{13.2-eq2},
we obtain
\begin{eqnarray}\label{13.2-eq5}
& & \3n\3n\3n\3n\3n\uplim_{\e\to 0}\mathbb{E}_{\omega_0}\frac{1}{\e}
\big(V(t_0+\e,\overline X(t_0+\e))-V(t_0,\overline X(t_0))\big)
\\
& & \3n\3n\3n\3n\3n\leq\!   \overline{R}(t_0,\omega_0)\!+ \!\lan
\bar{p}(t_0,\omega_0),A \overline X(t_0)\ran_H\! +\!\lan
\bar{p}(t_0,\omega_0),\bar{a}(t_0,\omega_0)\ran_H\!+\!\frac{1}{2}\lan
\overline{P}(t_0,\omega_0)\bar{b}(t_0,\omega_0),\bar{b}(t_0,\omega_0)
\ran_{\mathcal{L}_2^0}.\nonumber
\end{eqnarray}
Using \eqref{13.2-eq5} and {\bf Claim 1} again, by  Fatou's lemma, we get
\begin{eqnarray}\label{13.2-eq6}
& & \uplim_{\e\to
0^+}\mathbb{E}\frac{1}{\e}\left(V\big(t_0+\e,\overline
X(t_0+\e)\big)-V\big(t_0,\overline
X(t_0)\big)\right)\nonumber
\\
& & = \uplim_{\e\to 0^+}\mathbb{E}
\Big[\mathbb{E}_{\omega_0}\frac{1}{\e}\big(V\big(t_0+\e,\overline
X(t_0+\e)\big)-V\big(t_0,\overline
X(t_0)\big)\big)\Big]
\nonumber\\
& & \leq \mathbb{E}\Big\{\uplim_{\e\to
0^+}\Big[\mathbb{E}_{\omega_0}
\frac{1}{\e}\big(V\big(t_0+\e,\overline
X(t_0+\e)\big)-V\big(t_0,\overline
X(t_0)\big)\big)\Big]\Big\}
\\
& & \leq\mathbb{E}\Big(\overline{R}(t_0)+\lan
\bar{p}(t_0),A\overline X(t_0)\ran_H+\lan
\bar{p}(t_0),\bar{a}(t_0)\ran_H+\frac{1}{2}\lan
\overline{P}(t_0)\bar{b}(t_0),\bar{b}(t_0)\ran_{\mathcal{L}_2^0}\Big)
\nonumber
\end{eqnarray}
for a.e. $t_0\in [0,T-\d)$. Applying Lemma
\ref{lm7.1} to  $g(t)=\mathbb{E}V(t,\overline
X(t))$, and using \eqref{13.2-eq6} and
\eqref{13.2-eq0}, we obtain
\begin{equation}\label{12.15-eq7}
\mathbb{E}V\big(T-\d,\overline
X(T-\d)\big)-V(0,\eta)\leq
-\mathbb{E}\int_0^{T-\d} \bar{f}(t)dt.
\end{equation}
Noting that $V(\cd,\cd)$ is continuous and
$V(T,\overline X(T))=h(\overline{X}(T))$,
letting $\d\to 0$ in \eqref{12.15-eq7}, we
obtain that
$$
\mathbb{E}\Big(\int_0^T f\big(t,\overline
X(t),\bar{u}(t)\big)dt+h\big(\overline
X(T)\big)\Big)\leq V(0,\eta),
$$
which means that the control $\bar{u}(\cdot)$ is optimal.
\endpf

\ms

{\it Proof of Proposition \ref{cor3.2}}.
Obviously, \eqref{13.2-eq07} implies
\eqref{13.2-eq0}. We only need to prove that
\eqref{13.2-eq0} implies \eqref{13.2-eq07}.

Suppose \eqref{13.2-eq07} holds. By Proposition
\ref{prop-7.1}, for a.e.  $(t,\omega)\in
[0,T]\times H$  such that
$(\overline{R}(t,\omega), \bar{p}(t,\omega),$
$\overline R(t,\omega)) \in
D_{t+,x}^{1,2,+}V(t,\overline X(t)),\ \overline
X(t,\omega)=x$, there exists a function $\varphi
\in C^{1,2}([0,T]\times H)$ so that
\begin{equation}\label{13.2-eq08}
\begin{cases}\ds
\big(\varphi(t,x),\varphi_t(t,x),\varphi_x(t,x),\varphi_{xx}(t,x)\big)=\big(V(t,x),\overline R(t,\omega),\bar p(t,\omega),\overline P(t,\omega)\big),\\
\ns\ds \varphi(s,y)>v(t,x), \qquad \forall \
(t,x)\neq (s,y)\in [t,T]\times H.
\end{cases}
\end{equation}
Fix a $u\in U$. Let $X(\cdot)=X(\cdot;t,x,u)$ be the trajectory with the control $u(r)\equiv u$. Then by It\^o's formula, for any $s>t$ with $s-t>0$ small enough, we have
\begin{eqnarray}\label{12.15-eq9}
0 \3n& \leq \3n& \frac{1}{s-t}\mathbb{E}\big(V(t,x)-\varphi
(t,x)-V(s,X(s))+\varphi(s,X(s))\big) \nonumber
\\
\3n& \leq \3n&\frac{1}{s-t}\Big(\mathbb{E}\int_t^s f(r,X(r),u)dr-\varphi(t,x)+\varphi (s,X(s))\Big)\nonumber\\
\3n& = \3n&
\frac{1}{s-t}\mathbb{E}\int_t^s
\Big(\varphi_t(r,X(r))+\langle
\varphi_x(r,X(r)),AX(r)\rangle_H
\nonumber\\
& & \hspace{2.2cm}+ G(r,X(r),u_{\e,s}(r),\varphi_x(r,X(r)),\varphi_{xx}(r,X(r)))\Big)dr\nonumber\\
\end{eqnarray}
This leads to
$$\varphi_t(t,x)+\langle \varphi_x(t,x),Ax\rangle_H+G(t,x,u,\varphi_x(t,x),\varphi_{xx}(t,x))\ge 0,\ \ \forall u\in U.$$
Hence
$$\varphi_t(t,x)+\langle \varphi_x(t,x),Ax\rangle_H+\inf\limits_{u\in U}G(t,x,u,\varphi_x(t,x),\varphi_{xx}(t,x))\ge 0.$$
This, together with \eqref{13.2-eq08}, implies
that
$$
\begin{array}{ll}\ds
\varphi_t(t,\overline X(t,\omega)) +  \lan
\varphi_x(t,\overline X(t,\omega)), A\overline
X(t,\omega)\ran_H\\
\ns\ds + \inf\limits_{u\in U}G(t,\overline
X(t,\omega),u,\varphi_x(t,\overline
X(t,\omega)),\varphi_{xx}(t,\overline
X(t,\omega)))\ge 0,
\end{array}
$$
which yields
$$\overline R(t,\omega)+\lan \bar{p}(t,\omega),
AX(t,\omega)\ran_H+\inf\limits_{u\in
U}G(t,\overline X(t,\omega),u,\bar
p(t,\omega),\overline P (t,\omega))\ge 0$$
This combining with \eqref{13.2-eq0} gives
\eqref{13.2-eq07}.
\endpf

\section{Lipschitz continuity of the value function}\label{sec-lipschitz}

In this section, we are to prove Theorem
\ref{th4.1}. We first introduce an auxiliary
control problem to be used in the sequel. .

\subsection{An auxiliary control problem}\label{subsec-weakform}

Recall that for \textbf{Problem}
$\boldsymbol{(S_{t\eta})}$),  both the
probabililty space $(\Omega, \mathcal{F},
\mathbb{P})$ and the Brownian motion $W(\cd)$ on
$(\Omega, \mathcal{F}, \mathbb{P})$ are given a
priori, and our controls are  $\mathbf
F$-adapted processes. In this subsection, we
introduce a family of auxiliary control problems
in which only the probability space $(\Omega,
\mathcal{F}, \mathbb{P})$ is fixed and the
Brownian motion is  part of the controls. We
will see this newly introduced admissible
control is closely related to our original one
under the original strong formulation and plays
an important role in the proof of the Lipschitz
continuity of the value function.

Let $t\in [0,T)$,  denote by $\widetilde{\dbW}_t$ the set of all
cylindrical Brownian motions $\widetilde W_t(\cd)$ on $(\Omega,
\mathcal{F}, \mathbb{P})$ over $[t,T]$  (with $\widetilde W_t(t)=0$
almost surely). It is well known that $\widetilde W_t(\cd)$ is
continuous, $\dbP$-a.s. and admits a modification which is
continuous for all $\om\in\Om$. In what follows, we always take the
continuous modification of cylindrical Brownian motion.

For a given $\widetilde W_t(\cd)\in \widetilde{\dbW}_t$, write
$\textbf{F}_{\widetilde W_t}$ for the natural filtration generated
by $\widetilde W_t(\cd)$. Let
$$
\widetilde{\mathcal{U}}_{\widetilde
W_t}[t,T]\deq \left\{u:[t,T]\times\Omega\to
U\big| u \mbox{ is $\textbf{F}_{\widetilde
W_t}$-adapted}\right\}.
$$
Clearly, both $\textbf{F}_{\widetilde W_t}$ and
$\widetilde{\mathcal{U}}_{\widetilde W}[t,T]$ depend on the Brownian
motion $\widetilde W_t(\cd)$.

The admissible control set is
$$
\widetilde{\mathcal{U}}_{EX}[t,T]\deq
\left\{(\tilde u(\cdot),\widetilde W_t(\cdot)):
\widetilde W_t(\cdot)\in \widetilde{\dbW}_t,\;
\tilde u(\cdot)\in
\widetilde{\mathcal{U}}_{\widetilde
W_t}[t,T]\right\}.
$$

Consider the following control system:
\begin{equation}\label{4.1-eq6}
\begin{cases}
\ds d\widetilde X(s)=\big(A\widetilde X(s)+a(s,\widetilde
X(s),\tilde u(s))\big)ds
+b(s,\widetilde X(s),\tilde u(s))d\widetilde W_t(s), &s \in (t,T],\\
\ns\ds  \widetilde X(t)=\eta,
\end{cases}
\end{equation}
where $\eta\in H$,  and $(\tilde u(\cdot),\widetilde W(\cdot))\in
\widetilde{\mathcal{U}}_{EX}[t,T]$.

\ss

Under Assumption {\bf (S1)}, for any $\eta\in H$, \eqref{4.1-eq6}
admits a unique mild solution $\widetilde X(\cdot)$ (e.g.,
\cite[Theorem 3.14]{Lu2021}). Then for any $(t,\eta)\in [0,T]\times
H$ and $(\tilde u(\cdot),\widetilde W(\cdot))\in
\widetilde{\mathcal{U}}_{EX}[t,T]$, the cost functional
\begin{equation}
\widetilde{\mathcal{J}}(t,\eta;\tilde
u(\cdot))=\mathbb{E}\Big(\int_t^T f(s,\widetilde
X(s),\tilde u(s))ds+h(\widetilde X(T))\Big)
\end{equation}
is well-defined. So does the corresponding value function
\begin{equation}
\widetilde V(t,\eta)\deq \inf\limits_{(\tilde
u(\cdot),\widetilde W(\cdot))\in
\widetilde{\mathcal{U}}_{EX}[t,T]}
\widetilde{\mathcal{J}}\big(t,\eta;\tilde
u(\cdot)\big),\hspace{1cm} \forall\ (t,\eta)\in
[0,T]\times H.
\end{equation}

\ss

\begin{remark}\label{rmk1}
Compared with \textbf{Problem}
$\boldsymbol{(S_{t\eta})}$, we enlarge the
admissible control set
$\widetilde{\mathcal{U}}_{EX}[t,T]$ to admit the
Brownian as part of the control. But the
probability space  $(\Omega, \mathcal{F},
\mathbb{P})$ is fixed. Recall that there is
another formulation (which is called {\it weak}
formulation) in the literature, in which  the
probability space  $(\Omega, \mathcal{F},
\mathbb{P})$ is also part of the control (e.g.
\cite[Subsection 2.1.2]{Fabbri2017} or
\cite[Subsection 4.2, Chapter 2]{Yong1999}). For
the {\it weak} formulation, one can also define
the value function $V_w(\cd,\cd)$ (see
\cite[Subsection 2.1.2]{Fabbri2017} for the
details).
\end{remark}
%

%
%

By the definition of $V(\cd,\cd)$ and
$\widetilde V(\cd,\cd)$, it is clear that
\begin{equation}\label{4.1-eq8.1}
V(t,\eta)\geq \widetilde V(t,\eta),\ \ \ \forall \ (t,\eta)\in
[0,T]\times H.
\end{equation}
On the other hand, by  \cite[Theorem 2.22,
Chapter 2]{Fabbri2017},  for all $(t,\eta)\in
[0,T]\times H$,  $V(t,\eta)$ equals the value
function under the {\it  weak} formulation
mentioned in Remark \ref{rmk1}. Consequently, we
have
\begin{equation}\label{4.1-eq8}
V(t,\eta)= \widetilde V(t,\eta),\ \ \ \forall \ (t,\eta)\in
[0,T]\times H.
\end{equation}

\ss

Next, we introduce  a special case of the auxiliary control problem
presented above, which plays a major role in the next subsection.

For each $t>0$,  denote by $\mathcal{F}^t_r$ the
$\si$-algebra generated by
$\{W(\tau)-W(t)\}_{t\leq \tau\leq r}$ and by
$\mathbf F^t\deq\{\mathcal{F}^t_r\}_{t\leq r\leq
T}$ the natural filtration of the Brownian
motion $\{W(r)-W(t)\}_{t\leq r\leq T}$. Write
$\mathbb{F}^t$ for the progressive
$\sigma$-algebra with respect to $\textbf{F}^t$.
Let
$$
\mathcal{U}^t[t,T] \deq
\left\{u(\cdot)\in\mathcal{U}[t,T] \big| \
u(r)\text{ is $\mathbf F^t$-adapted, } \forall\
t\leq r\leq T\right\},
$$

For $u(\cdot)\in \mathcal{U}^t[t,T]$ and any $s\in [0,T)$, let
$$
\tilde{u}(\cdot)=u(\tau(\cdot)), \q
\widetilde{W}(\cdot)=\sqrt{1/\dot{\tau}}W(\tau(\cdot))-W(t).
$$
where $\tau(r)=\frac{T(t-s)+(T-t)r}{T-s}$.

We claim that $(\tilde{u}(\cdot),
\widetilde{W}(\cdot))\in
\widetilde{\mathcal{U}}_{EX}[s,T]$. Indeed,  it
is clear that

(i) $\widetilde
W(s)=\sqrt{1/\dot{\tau}}W(t)-W(t)=0$,  for that
$W(\cdot)-W(t)$ is cylindrical Brownian motion;

(ii) For all $n\in \mathbb{N}^+$ and  $s=s_1<s_2<\cdots<s_n\leq T$,
we have
\begin{eqnarray*}
& & \left( \widetilde W(s_1),\ \widetilde
W(s_2)-\widetilde W(s_1),\cdots,
\widetilde W(s_n)-\widetilde W(s_{n-1})\right)\\
\3n& = \3n& \left(
\sqrt{1/\dot{\tau}}W(r_1))-W(t),
\sqrt{1/\dot{\tau}}W(r_2)-\sqrt{1/\dot{\tau}}W(r_1),\cdots,
\sqrt{1/\dot{\tau}}W(r_n)-\sqrt{1/\dot{\tau}}W(r_{n-1})
\right)
\end{eqnarray*}
are independent, where $r_i=\tau(s_i),\ \forall \ 1\leq i\leq n$.

(iii) For any $s\leq r<l\leq T$,
$$\widetilde W(r)-\widetilde W(l)=\sqrt{1/\dot{\tau}}\big(W(\tau(r))-W(\tau(l))\big)
\sim N(0, (1/\dot{\tau}(\tau(r)-\tau(l)))I)=N(0,(r-l)I),
$$
where $I$ is the identity operator in $\wt H$. On the other hand,
\begin{eqnarray}\label{12.15-eq1}
\widetilde{\mathcal{F}}_{s}^r \3n&\triangleq \3n&
\sigma\big\{\widetilde{W}(l):\ s\leq l\leq r\big\} =
\sigma\big\{\widetilde{W}^{-1}(B):\ B\in \mathcal{B}(\mathbb{R}),\
s\leq l\leq r\big\}\nonumber
\\
\3n&=\3n& \sigma\big\{\big(\sqrt{\dot{\tau}}W^{-1}(\tau(l))-W(t)\big)(B): \ B\in \mathcal{B}(\mathbb{R}),\ s\leq l\leq r\big\}\\
\3n&=\3n& \sigma \big\{\big(W^{-1}(\rho)-W(t)\big)(B):\ B\in
\mathcal{B}(\mathbb{R}),\ t\leq \rho\leq \tau(r)\big\}=
\mathcal{F}^{\tau(r)}_t.\nonumber
\end{eqnarray}
Since $\tilde u(r)=u(\tau(r))$ is $\mathcal{F}^{\tau(r)}_t$-
measurable, it follows from \eqref{12.15-eq1} that $\tilde
u(r)=u(\tau(r))$ is $\widetilde{\mathcal{F}}_{s}^r$- measurable.

We  end up this subsection with the following result.
\begin{proposition}\label{prop3.1-4}
For any $\eta\in H,\ t\in [0,T)$,
\begin{equation}\label{12.13-eq6}
\inf\limits_{u(\cdot)\in\mathcal{U}[t,T]}
{\mathcal{J}}(t,\eta;u(\cdot))=\inf\limits_{u(\cdot)\in\mathcal{U}^t[t,T]}
{\mathcal{J}}(t,\eta;u(\cdot)).
\end{equation}
\end{proposition}
The proof of Proposition \ref{prop3.1-4} is very similar to the one
for \cite[Subsection 4.1, Proposition 4.3]{Chen2022}. We provide it
here for the convenience of readers.

\ss

{\it Proof}. We divide the proof into three steps.

\textbf{Step\ 1.} Let
$$
\mathcal{U}_D^t \deq \Big\{u(s)=\sum\limits_{j=1}^N
u^j(s)1_{\Omega_j}\Big|\, u^j(s)\in\mathcal{U}^t[t,T],\;
\{\Omega_j\}_{j=1}^N\subset \mathcal{F}_t\ \text {is a partition of
} \Omega\Big\}.
$$
Since $\mathcal{U}_D^t\subset \mathcal{U}[t,T]$, we have
\begin{equation}\label{12.13-eq4-1-1}
\operatorname*{inf}_{u(\cdot)\in\mathcal{U}[t,T]}
{\mathcal{J}}(t,\eta;u(\cdot))\leq \operatorname*{inf}_{u(\cdot)\in
\mathcal{U}_D^t} {\mathcal{J}}(t,\eta;u(\cdot)).
\end{equation}
On the other hand, from \cite[Lemma 4.12]{Yong2020}, we know that
$\mathcal{U}_D^t$ is dense in $\mathcal{U}[t,T]$. Consequently,
\begin{equation}\label{12.13-eq4-1}
\operatorname*{ess~inf}_{u(\cdot)\in\mathcal{U}_D^t}
{\mathcal{J}}(t,\eta;u(\cdot))\leq
\operatorname*{ess~inf}_{u(\cdot)\in\mathcal{U}[t,T]}
{\mathcal{J}}(t,\eta;u(\cdot)).
\end{equation}
From \eqref{12.13-eq4-1-1} and \eqref{12.13-eq4-1}, we see that
\begin{equation}\label{12.13-eq4}
\operatorname*{inf}_{u(\cdot)\in\mathcal{U}[t,T]}
{\mathcal{J}}(t,\eta;u(\cdot))=\operatorname*{inf}_{u(\cdot)\in
\mathcal{U}_D^t} {\mathcal{J}}(t,\eta;u(\cdot)).
\end{equation}

\ss

\textbf{Step\ 2.} In this step, we prove
\begin{equation}\label{12.13-eq5}
\operatorname*{inf}_{u(\cdot)\in\mathcal{U}_D^t}
{\mathcal{J}}(t,\eta;u(\cdot))=\operatorname*{inf}_{u(\cdot)\in\mathcal{U}^t[t,T]}
{\mathcal{J}}(t,\eta;u(\cdot)).
\end{equation}
Since $\mathcal{U}^t\subset \mathcal{U}_D^t$, we have
\begin{equation}\label{12.13-eq5-1}
\operatorname*{inf}_{u(\cdot)\in\mathcal{U}_D^t}
{\mathcal{J}}(t,\eta;u(\cdot))\leq
\operatorname*{inf}_{u(\cdot)\in\mathcal{U}^t[t,T]}
{\mathcal{J}}(t,\eta;u(\cdot)).
\end{equation}

Now we show the inverse inequality of \eqref{12.13-eq5-1}.
For all $u(\cdot)\in\mathcal{U}_D^t$, we have
$$
{\mathcal{J}}(t,\eta;u(\cdot))=
{\mathcal{J}}\Big(t,\eta;\sum\limits_{j=1}^N\chi_{\Omega_j}u^j(\cdot)\Big)
=\sum\limits_{j=1}^N\chi_{\Omega_j}
{\mathcal{J}}(t,\eta;u^j(\cdot)).
$$
For $j=1,2,\cdots,N$, noting that $u^j(\cdot)$ is $\dbF^t$
measurable,  we find that ${\mathcal{J}}(t,\eta;u^j(\cdot))$ is
deterministic. Without loss of generality, we assume that
$$
{\mathcal{J}}(t,\eta;u^1(\cdot))\leq
{\mathcal{J}}(t,\eta;u^j(\cdot)),\quad \forall j=2,3,\cdots,N.
$$
Thus, it holds that
$$
{\mathcal{J}}(t,\eta;u(\cdot))\geq
{\mathcal{J}}(t,\eta;u^1(\cdot))\geq
\operatorname*{inf}_{u(\cdot)\in\mathcal{U}^t[t,T]}
{\mathcal{J}}(t,\eta;u(\cdot)).
$$
For $u(\cdot)\in\mathcal{U}^t[t,T]$ is
arbitrarily chosen, we get that
$$
\operatorname*{inf}_{u(\cdot)\in\mathcal{U}_D^t}
{\mathcal{J}}(t,\eta;u(\cdot))\geq
\operatorname*{inf}_{u(\cdot)\in\mathcal{U}^t[t,T]}
{\mathcal{J}}(t,\eta;u(\cdot)).
$$
\ss

\textbf{Step\ 3.} We finish the proof in this step.

From \eqref{12.13-eq4} and \eqref{12.13-eq5}, we see that
\begin{equation}\label{12.13-eq6-1}
\operatorname*{inf}\limits_{u(\cdot)\in\mathcal{U}[t,T]}
{\mathcal{J}}(t,\eta;u(\cdot))=\operatorname*{inf}\limits_{u(\cdot)\in\mathcal{U}^t[t,T]}
{\mathcal{J}}(t,\eta;u(\cdot)).
\end{equation}
The right hand side of \eqref{12.13-eq6-1} is deterministic, hence
\eqref{12.13-eq6-1} can be simplified as
$$
\inf\limits_{u(\cdot)\in\mathcal{U}[t,T] }
{\mathcal{J}}(t,\eta;u(\cdot))
=\inf\limits_{u(\cdot)\in\mathcal{U}^t[t,T]}
{\mathcal{J}}(t,\eta;u(\cdot)).
$$
\endpf

\subsection{Proof of Theorem \ref{th4.1}}

Before completing the proof of Theorem
\ref{th4.1}, we present some preliminaries to be
used latter, which will play a crucial role in
the proof of the regularity result. More details
can be found in \cite{Pazy1983}.

\ss

\begin{proposition}\label{prop1.1}
Suppose $0<t_1<t_2,\ \eta\in H$. Then
\begin{equation*}
\big|S(t_2)\eta-S(t_1)\eta\big|_H\leq
(t_2-t_1)\big|AS(t_1)\big|_{\cL(H)}|\eta|_H.
\end{equation*}
\end{proposition}

Proposition \ref{prop1.1} should be a well-known result and its
proof is very easy. However, we do not find an exact reference for
it. For the convenience of readers, we provide a proof below.

\ss

{\it Proof}. Since $\{S(t)\}_{t\ge 0}$ is analytic and contractive,
the map $t\to S(t)\eta,\ \forall \ \eta\in H$, is differentiable for
$t>0$ and $\big|S(t)\big|\leq 1$. Then $AS(t_1)$ is a bounded linear
operator and
\begin{eqnarray*}
\big|S(t_2)\eta-S(t_1)\eta\big|_H\3n & = &
\3n\Big|\int_{t_1}^{t_2} AS(t)\eta dt\Big|_H
=\Big|\int_{t_1}^{t_2} S(t-t_1)AS(t_1)\eta dt \Big|_H\\
& \leq & (t_2-t_1)\big|AS(t_1)\big|_{\cL(H)}|\eta|_H.
\end{eqnarray*}
\ss

The next proposition is taken from \cite[Theorem 6.13, Chapter 2]{Pazy1983}.

\ss

\begin{proposition}\label{prop1.2}
There exists a constant $C>0$ such that
\begin{equation*}
\big|A S(t)\big|_{\cL(H)}\leq \frac{C}{t},
\hspace{1cm} \forall\ t>0.
\end{equation*}
\end{proposition}

Now let us prove Theorem \ref{th4.1}.

\ss

{\it Proof of Theorem \ref{th4.1}}. Let us fix $\delta>0$, and let
$(t_1,\eta_1)$ and $(t_0,\eta_0)\in [0,T)\times H$ be such that
\begin{equation}\label{4-eq1}
\min \{T-t_1,T-t_0\}>\delta.
\end{equation}
Without loss of generality,  assume that
$t_0>t_1$.

For any $\varepsilon>0$, by Proposition \ref{prop3.1-4}, there
exists $u_0(\cdot)\in \mathcal{U}^{t_0}[t_0,T]$ such that
\begin{equation}\label{12.15-eq2}
\mathcal{J}(t_0,\eta_0;u_0(\cdot)) < V(t_0,\eta_0) + \varepsilon.
\end{equation}
Let $X_0(\cdot)$ be the solution of
\begin{equation}\label{4-eq3}
\begin{cases}
\ds dX_0(t)=\big(AX_0(t)+a(t,X_0(t),u_0(t))\big)dt+b(t,X_0(t),u_0(t))dW(t), &t \in (t_0,T],\\
\ns\ds   X(t_0)=\eta_0.
\end{cases}
\end{equation}
Consider the change of time as follows:
\begin{equation}\label{4-eq4}
\tau:[t_1,T]\to [t_0,T],\ \
\tau(t)=\frac{T(t_0-t_1)+(T-t_0)t}{T-t_1}.
\end{equation}
Then it holds that
\begin{equation}\label{4-eq5}
\dot{\tau}(t)=\frac{T-t_0}{T-t_1},\ \
\tau^{-1}(t)=\frac{t(T-t_1)-T(t_0-t_1)}{T-t_0},
\end{equation}
and that
\begin{equation}\label{4-eq6}
\tau(t_1)=t_0,\ \ \tau(T)=T.
\end{equation}
Let $\tilde u(t)=u_0(\tau(t))$, and denote by $\widetilde X(t)$ the
solution of
\begin{equation}\label{4-eq7}
\begin{cases}
d\widetilde X(t)=\big(A\widetilde
X(t)+a\big(t,\widetilde X(t),\tilde
u(t)\big)\big)dt+ b\big(t,\widetilde X(t),\tilde
u(t)\big)d\widetilde W(t),
&t \in (t_1,T],\\
\ns\ds \wt X(t_1)=\eta_1,
\end{cases}
\end{equation}
where $\widetilde W(t)=\sqrt{1/\dot{\tau}}W(\tau(t)).$ Obviously,
$(\widetilde W(\cdot),\tilde u(\cdot))\in
\widetilde{\mathcal{U}}_{EX}[t_1,T]$. Thus,
\begin{eqnarray}\label{4-eq8}
\wt X(t) \3n& = &\3n S(t-t_1)\eta_1 + \int_{t_1}^t
S(t-s)a\big(s,\wt X(s),\tilde u(s)\big)ds
+\int_{t_1}^t S(t-s)b\big(s,\wt X(s),\tilde u(s)\big)d\wt W(s)\nonumber\\
& = &\3n S(t-t_1)\eta_1 + \int_{t_1}^t S(t-s)
a\big(s,\wt X(s),u_0(\tau(s))\big)ds\\
& & \q+\int_{t_1}^t \sqrt{1/\dot{\tau}}
S(t-s)b\big(s,\wt
X(s),u_0(\tau(s))\big)dW(\tau(s)).\nonumber
\end{eqnarray}

Next, by  \eqref{4.1-eq8} and \eqref{12.15-eq2}, we have
\begin{eqnarray}\label{12.15-eq5}
& & V(t_1,\eta_1)-V(t_0,\eta_0)- \varepsilon
=\widetilde V(t_1,\eta_1)-V(t_0,\eta_0)- \varepsilon \nonumber \\
& & \leq \widetilde{\mathcal{J}}(t_1,\eta_1;\tilde
u)-\mathcal{J}(t_0,\eta_0;u_0)
\\
& & = \mathbb{E}\Big(\int_{t_1}^T f(t,\widetilde X(t),\tilde
u(t))dt-\int_{t_0}^T f(t,X_0(t),u_0(t))dt + h(\wt
X(T))-h(X_0(T))\Big). \nonumber
\end{eqnarray}

Set $X_1(r)=\widetilde X(\tau^{-1}(r))$ for $r\in [t_1,T]$. By
changing the variable $r=\tau(t)$ in the first  integral of
\eqref{12.15-eq5}, we obtain that
\begin{eqnarray*}
& & V(t_1,\eta_1)-V(t_0,\eta_0)-\varepsilon\\
& & \leq \mathbb{E} \int_{t_0}^T
\Big(\frac{1}{\dot{\tau}}
f\big(\tau^{-1}(r),X_1(r),u_0(r)\big)-
f\big(r,X_0(r),u_0(r)\big)\Big)dr +
\mathbb{E}\big( h(X_1(T))-h(X_0(T))\big).
\end{eqnarray*}
By Assumption {\bf{(S2)$'$}},  we have
\begin{eqnarray}\label{12.15-eq3}
& &\3n\3n V(t_1,\eta_1)-V(t_0,\eta_0)-\varepsilon \nonumber\\
& & \3n\3n\leq \mathbb{E} \int_{t_0}^T \Big(\frac{1}{\dot{\tau}}
f(\tau^{-1}(r),X_1(r), u_0(r))-f(\tau^{-1}(r),X_1(r),u_0(r))\Big)dr
\\
& &  +\mathbb{E} \int_{t_0}^T
\big(f(\tau^{-1}(r),X_1(r),u_0(r))-f(r,X_0(r),u_0(r))\big)dr +
\mathbb{E}\big( h(X_1(T))-h(X_0(T))\big)  \nonumber\\
& & \3n\3n\leq C\mathbb{E} \int_{t_0}^T
\Big(\Big|1-\frac{1}{\dot{\tau}}\Big|
+\big|\tau^{-1}(r)-r\big|+\big|X_1(r)-X_0(r)\big|_H\Big)dr +
C\mathbb{E}\big|X_1(T)-X_0(T)\big|_H. \nonumber
\end{eqnarray}

We claim the following estimates hold:
\begin{eqnarray}
\Big|1-\frac{1}{\dot{\tau}}\Big| & \leq & C|t_1-t_0|,\label{4-eq9}\\
|\tau^{-1}(r)-r| & \leq & C|t_1-t_0|\label{4-eq10},\\
\mathbb{E} \int_{t_0}^T \big|X_1(r)-X_0(r)\big|_Hdr &\leq &
C\big(|t_1-t_0|+|\eta_1-\eta_0|_H\big)\big(1+|\eta_1|_H\big),
\label{4-eq11}\\
\mathbb{E}\big|X_1(T)-X_0(T)\big|_H & \leq &
C\big(|t_1-t_0|+|\eta_1-\eta_0|_H\big)\big(1+|\eta_1|_H\big).\label{4-eq12}
\end{eqnarray}

{\it Proof of \eqref{4-eq9}}. By the definition
of  $\tau(\cd)$,
$$
\Big|1-\frac{1}{\dot{\tau}}\Big|
=\Big|1-\frac{T-t_1}{T-t_0}\Big|=\Big|\frac{t_1-t_0}{T-t_0}\Big|\leq
\frac{1}{\delta}|t_1-t_0|\leq C|t_1-t_0|.
$$

{\it Proof of \eqref{4-eq10}}. By the definition
of $\tau(\cd)$ and \eqref{4-eq1}, we obtain
$$|\tau^{-1}(r)-r|=\Big|\frac{(T-t_1)r
-(t_0-t_1)T}{T-t_0}-r\Big|=\Big|\frac{(t_0-t_1)(T-r)}{T-t_0}\Big|\leq
C|t_1-t_0|.
$$

{\it Proof of \eqref{4-eq11}}. Recalling   the
definition of $X_1(r)$, we conclude that
\begin{eqnarray*}
& & \3n\3n \mathbb{E}\big|X_1(r)-X_0(r)\big|_H^2
= \mathbb{E}\big|\widetilde X(\tau^{-1}(r))-X_0(r)\big|_H^2\\
& &\3n\3n =
\mathbb{E}\Big[\Big|S(\tau^{-1}(r)-t_1) \eta_1
-S(r-t_0)\eta_0 +\int_{t_1}^{\tau^{-1}(r)}
S(\tau^{-1}(r)-\rho)a\big(\rho,\widetilde
X(\rho),u_0(\tau(\rho))\big)d\rho
\\
& & \hspace{2mm}-\!\int_{t_0}^r\!\!
S(r\!-\!\rho)a\big(\rho,X_0(\rho),u_0(\rho)\big)d\rho\!
+\!\int_{t_1}^{\tau^{-1}(r)}\!\!\!\!
\sqrt{\frac{1}{\dot{\tau}}}
S(\tau^{-1}\!(r)\!-\!\rho)b\big(\rho,\widetilde
X(\rho),u_0 (\tau(\rho))\big)dW\!(\tau(\rho))
\\
& & \hspace{2mm}-\int_{t_0}^r S(r-\rho) b\big(\rho,X_0(\rho),u_0(\rho)\big)dW(\rho)\Big|_H^2\Big]\\
& &\3n\3n =\mathbb{E}\Big\{\Big|S(\tau^{-1}(r)-t_1)\eta_1-S(r-t_0)\eta_0\\
& & \hspace{2mm} +
\int_{t_0}^r\Big[\frac{1}{\dot{\tau}}S(\tau^{-1}(r)\!
-\!\tau^{-1}(\rho))a(\tau^{-1}(\rho), X_1(\rho),u_0(\rho))-
S(r\!-\!\rho)a(\rho,X_0(\rho),u_0(\rho))\Big]d\rho
\\
& & \hspace{2mm} +
\int_{t_0}^r\big[\dot{\tau}^{-1/2}S(\tau^{-1}(r)\!-\!
\tau^{-1}(\rho))b\big(\tau^{-1}(\rho),
X_1(\rho),u_0(\rho)\big)
\\
& & \hspace{1cm}
-S(r-\rho)b\big(\rho,X_0(\rho),u_0(\rho)\big)\big]dW(\rho)\Big|_H^2\Big\}
\end{eqnarray*}
By Assumption  {\bf{(S1)$'$}},  we have
\begin{eqnarray*}
& & \mathbb{E}\big|X_1(r)-X_0(r)\big|_H^2 \\
& & \leq C\mathbb{E}\big(\big|S(\tau^{-1}(r)-t_1)\eta_1
-S(r-t_0)\eta_1\big|_H^2+\big|S(r-t_0)\big|_{\cL(H)}^2\big|\eta_1
-\eta_0\big|_H^2\big)
\\
& & \q + C\mathbb{E}
\Big[\int_{t_0}^r\big(\big|\dot{\tau}^{-1}-1\big|
\big|S(\tau^{-1}(r)-\tau^{-1}(\rho))\big|_{\cL(H)}
\big|a(\tau^{-1}(\rho),X_1(\rho),u_0(\rho))\big|_H
\\
& & \hspace{1.3cm} +
\big|S(\tau^{-1}(r)-\tau^{-1}(\rho))\big|_{\cL(H)}
\big|a(\tau^{-1}(\rho),X_1(\rho),u_0(\rho))-a(\rho,X_0(\rho),u_0(\rho))\big|_H
\\
& & \hspace{1.3cm} + \big|\big[S(\tau^{-1}(r)-\tau^{-1}(\rho))
-S(r-\rho)\big] a(\rho,X_0(\rho),u_0(\rho))\big|_H\big)d\rho\Big]^2
\\
& & \q + C\mathbb{E} \int_{t_0}^r
\big(\big|\dot{\tau}^{-1/2}-1\big|\big|S(\tau^{-1}(r)-
\tau^{-1}(\rho))\big|_{\cL(H)}\big|b(\tau^{-1}(\rho),X_1(\rho),u_0(\rho))\big|_H
\\
& & \hspace{1.3cm} +
\big|S(\tau^{-1}(r)-\tau^{-1}(\rho))\big|_{\cL(H)}
\big|b(\tau^{-1}(\rho),X_1(\rho),u_0(\rho))-b(\rho,X_0(\rho),u_0(\rho))\big|_H
\\
& & \hspace{1.3cm} +
\big|(S(\tau^{-1}(r)-\tau^{-1}(\rho))-S(r-\rho))b(\rho,X_0(\rho),u_0(\rho))\big|_H\big)^2
d\rho
\\
& & \leq
C\Big(\Big|S\Big((r-t_0)\frac{T-t_1}{T-t_0}\Big)\eta_1-S(r-t_0)\eta_1\Big|_{\cL(H)}^2
+|\eta_1-\eta_0|_H^2\Big)
\\
& & \hspace{1cm} + C\mathbb{E}\Big(|t_1-t_0|^2+\int_{t_0}^r
|X_1(\rho)-X_0(\rho)|_H^2d\rho\Big)
\\
& & \hspace{1cm} + C\mathbb{E} \int_{t_0}^r \Big|
\Big[S\Big((r-\rho)\frac{T-t_1}{T-t_0}\Big)-
S(r-\rho)\Big]a(\rho,X_0(\rho),u_0(\rho))\Big|_H^2 d\rho
\\
& & \hspace{1cm} + C\mathbb{E} \int_{t_0}^r
\Big|\Big[S\Big((r-\rho)\frac{T-t_1}{T-t_0}\Big)-
S(r-\rho)\Big]b(\rho,X_0(\rho),u_0(\rho))\Big|_H^2 d\rho.
\end{eqnarray*}

Recalling that $t_0>t_1$, we have
$(r-\rho)\frac{T-t_1}{T-t_0}>r-\rho$. Noting that $\{S(t)\}_{t\ge
0}$ is an analytic semigroup, for any $\eta\in H$, we reach that
\begin{eqnarray}\label{12.15-eq4}
& & \Big|S\Big((r-\rho)\frac{T-t_1}{T-t_0}\Big)\eta -
S(r-\rho)\eta\Big|_H \nonumber\\
& & \leq
C(r-\rho)\frac{t_0-t_1}{T-t_0}\big|AS(r-\rho)\big|_{\cL(H)}|\eta|_H
\\
& & \leq C(r-\rho)\frac{t_0-t_1}{T-t_0}\frac{1}{r-\rho}|\eta|_H\leq
C(t_0-t_1)|\eta|_H.\nonumber
\end{eqnarray}
From \eqref{12.15-eq4}, we see that
\begin{eqnarray*}
& & \mathbb{E}\big|X_1(r)-X_0(r)\big|_H^2\\
& & \leq C\big(|\eta_1-\eta_0|_H^2 +|t_1-t_0|^2 |\eta_1|_H^2\big) +
C\mathbb{E}\Big(|t_1-t_0|^2+\int_{t_0}^r
|X_1(\rho)-X_0(\rho)|_H^2d\rho\Big)
\\
& & \leq
C\big[|\eta_1-\eta_0|_H^2+|t_1-t_0|^2\big(|\eta_1|_H^2+1\big)\big]
+C\mathbb{E} \int_{t_0}^r |X_1(\rho)-X_0(\rho)|_H^2d\rho.
\end{eqnarray*}
This, along with Gronwall's inequality, implies that
\begin{equation}\label{4-eq13}
\mathbb{E}\big|X_1(r)-X_0(r)\big|_H^2 \leq
C\big(|\eta_1-\eta_0|_H^2+|t_1-t_0|^2\big)
\big(|\eta_1|_H^2+1\big),
\end{equation}
which yields
\begin{equation}\label{4-eq14}
\mathbb{E}\big|X_1(r)-X_0(r)\big|_H\leq C\big(|\eta_1-\eta_0|_H
+|t_1-t_0|\big)\big(|\eta_1|_H+1\big).
\end{equation}

{\it Proof of \eqref{4-eq12}}. Owing to \eqref{4-eq11},
\begin{equation*}
\mathbb{E}\big|X_1(T)-X_0(T)\big|_H =\mathbb{E}\big|\widetilde
X(T)-X_0(T)\big|_H\leq C\big(|\eta_1-\eta_0|_H
+|t_1-t_0|\big)\big(|\eta_1|_H+1\big).
\end{equation*}

Combining \eqref{12.15-eq3}--\eqref{4-eq12}, we conclude  that
$$V(t_1,\eta_1)-V(t_0,\eta_0)-\e\leq C\big(|\eta_1|_H+1\big)\big(|t_1-t_0|+|\eta_1-\eta_0|_H\big).$$
From the arbitrariness of $\e$, we obtain that
$$V(t_1,\eta_1)-V(t_0,\eta_0) \leq C\big(|\eta_1|_H+1\big)\big(|t_1-t_0|+|\eta_1-\eta_0|_H\big).$$
Similarly, we can prove that
$$V(t_0,\eta_0)- V(t_1,\eta_1) \leq C\big(|\eta_1|_H+1\big)\big(|t_1-t_0|+|\eta_1-\eta_0|_H\big).$$
Therefore, we conclude that
$$|V(t_1,\eta_1)-V(t_0,\eta_0)|\leq C\big(|\eta_1|_H+1\big)\big(|t_1-t_0|+|\eta_1-\eta_0|_H\big).$$
This completes the proof.
\endpf

\begin{remark}
In the proof of Theorem \ref{th4.1}, we use the weakly formulated
admissible control $(\tilde u(\cdot), \tilde W(\cdot))$ and the
control system driven by it as an auxiliary tool to legitimate our
``change of time" strategy applied here. We borrow this idea from
\cite{Buckdahn2010}.
\end{remark}


\section{An Illustrative example }\label{sec-exam}

In this section, we present an illustrative
example which fulfill the assumptions in Theorem
\ref{th7.1}   and Corollary \ref{cor3.1}.

Let $\cO\subset\dbR^n$ be a bounded domain with the smooth boundary
$\pa \cO$.  Let $H=L^2(\cO)$ and $U$  be a bounded closed subset of
$L^2(\cO)$. Consider the following stochastic parabolic equation:
\begin{equation}\label{system3}
\begin{cases}
\ds   dy =\big(\Delta y + \tilde a(t,y,u)\big) dt+\tilde b(t,y,u) dW(t) &\textup{in } (0,T]\times \cO,\\
\ns\ds   y =0 &\textup{on }   (0,T]\times \pa \cO,\\
\ns\ds   y(0)=\eta &\textup{in }  \cO,
\end{cases}
\end{equation}
where $ \eta \in L^2(\cO)$, $u(\cdot)\in \cU[0,T]$, and $\tilde a$
and $\tilde b$  satisfy the following condition:

\ss

\no{\bf (B1)} {\it For $\f=\tilde a,\tilde b$, suppose that
    $\f(\cd,\cd,\cd):[0,T]\times \dbR\times \dbR\to
    \dbR$ satisfies : i) For any $(r,u)\in
    \dbR\times \dbR$, the function
    $\f(\cd,r,u):[0,T] \to \dbR$ is Lebesgue
    measurable; ii) For any $t\in [0,T],\ r\in
    \dbR$, the function $\f(t,r,\cd):\dbR\to \dbR$
    is continuous; and iii) For all
    $(t,t_1,t_2,r_1,r_2,u)\in [0,T]\times \dbR\times
    \dbR\times \dbR$,
    \begin{equation}\label{ab0}
        \left\{
        \begin{array}{ll}\ds
            |\f(t,r_1,u) - \f(t,r_2,u)|   \leq
            \cC (|t_1-t_2|+|r_1-r_2|,\\
            \ns\ds |\f(t,0,u)| \leq \cC;
        \end{array}
        \right.
\end{equation}
 iv)  For all
 $(t,u)\in [0,T]\times  \dbR$,
    $\f(t,\cd,u)$ are $C^2$, and for any $(r,u)\in
    \dbR\times \dbR$ and a.e. $t\in [0,T]$,
    \begin{equation*}\label{ab1}
            |\f_r(t,r,u)|
         \leq \cC.
\end{equation*}}

 Consider the following
cost functional:
\begin{equation}\label{cost3}
\mathcal{J}(\eta;u(\cdot))= \mathbb{E}\Big(\int_0^T\int_G \tilde
f(t,y(t),u(t))dxdt+\int_G \tilde h(y(T))dx \Big),
\end{equation}
where $\tilde f$ and $\tilde h$ satisfy the following condition:

\ss

\no{\bf (B2)} {\it $\tilde f(\cdot,r,u)$ is Lipschitz, $\tilde
f(t,\cd,u)$ and $\tilde h(\cd)$ are $C^2$, such that $\tilde
f_r(t,r,\cd)$ and $\tilde f_{rr}(t,r,\cd)$ are continuous, and for
any $(r,u)\in \dbR\times \dbR$ and a.e. $t\in [0,T]$,
\begin{equation*}\label{ab1-1}
|\tilde f_r(t,r,u)| +|\tilde h_r(r) | \leq \cC.
\end{equation*}
\ss

\no{\bf (B3)} {\it For $\f=\tilde a,\tilde b$, $\f(\cd,0,\cd)=0$.}
}

Under {\bf (B1)} and {\bf (B2)}, it is easy to see that {\bf
(S1)$'$}--{\bf (S2)$'$} hold. Under {\bf (B1)}--{\bf (B3)}, by the
regularity theory of stochastic parabolic equations
(e.g.,\cite{Flandoli1990}), we know that $\D y\in L^2_\dbF(0,T;H)$,
namely, \eqref{13.2-eq01} holds.  Hence,  all assumptions in Theorem
\ref{th7.1} and  Corollary \ref{cor3.1} are fulfilled.

\end{document}